\documentclass[11pt,a4paper,leqno,twoside]{amsart}
\usepackage{latexsym,amssymb,amsmath}
\input xy
\xyoption{all}

\def\trace{{\rm trace}}

\def\cE{{\mathcal E}}
\def\cF{{\mathcal F}}

\def\cO{{\mathcal O}}

\def\CC{\mathbb C}
\def\RR{\mathbb R}
\def\HH{\mathbb H}
\def\AA{{\mathbb A}}

\def\OO{\mathbb O}
\def\DD{\mathbb D}

\def\ZZ{\mathbb Z}

\def\11{\mathbf 1}
\def\PP{\mathbb P}
\def\QQ{\mathbb Q}
\def\FF{\mathbb F}

\def\e1{\varepsilon_1}
\def\e2{\varepsilon_2}
\def\e3{\varepsilon_3}

\def\P2{{\PP}^2}

\def\JA{{\mathcal J}_3(\AA)}

\def\JR{{\mathcal J}_3(\RR)}
\def\JC{{\mathcal J}_3(\CC)}
\def\JH{{\mathcal J}_3(\HH)}
\def\JO{{\mathcal J}_3(\OO)}
\def\00{\underline{0}}
\def\J0{{\mathcal J}_3(\underline{0})}

\def\PJA{\PP{\mathcal J}_3(\AA)}

\def\PJ0{\PP({\mathcal J}_3(\underline{0}))}

\def\fsl{{\mathfrak {sl}}}
\def\fgl{{\mathfrak {gl}}}
\def\fsp{{\mathfrak {sp}}}

\def\fso{{\mathfrak {so}}}
\def\fe{{\mathfrak e}}

\def\ff{{\mathfrak f}}

\def\ft{{\mathfrak t}}

\def\l{\lambda}
\def\a{\alpha}

\def\om{\omega}

\def\b{\beta}

\def\s{\sigma}

\def\th{\theta}
\def\m{\mu}

\def\e{\varepsilon}
\def\ot{{\mathord{\,\otimes }\,}}
\def\op{{\mathord{\,\oplus }\,}}

\def\lra{{\mathord{\;\longrightarrow\;}}}
\def\ra{{\mathord{\;\rightarrow\;}}}

\def\AP2{{\AA\PP}^2}
\def\RP2{{\RR\PP}^2}
\def\CP2{{\CC\PP}^2}
\def\HP2{{\HH\PP}^2}
\def\OP2{{\OO\PP}^2}
\def\DA{\DD_3({\AA})}
\def\DR{\DD_3({\RR})}
\def\DC{\DD_3({\CC})}
\def\DH{\DD_3({\HH})}

\def\SLA{SL_3({\AA})}

\def\SOA{SO_3({\AA})}
\def\fsoa{{\mathfrak so}_3(\AA)}

\def\tr{{\rm trace}\;}
\def\dim{{\rm dim}\;}

\newtheorem{theo}{Theorem}[section]
\newtheorem{coro}[theo]{Corollary}
\newtheorem{lemm}[theo]{Lemma}
\newtheorem{prop}[theo]{Proposition}
\newtheorem{defi}[theo]{Definition}

\begin{document}

\title{Severi varieties\linebreak and  their varieties of reductions}
\author{A. Iliev, L. Manivel}

\begin{abstract} We study the varieties of reductions associated to the four 
Severi varieties, the first example of which is the Fano threefold of index $2$ 
and degree $5$ studied by Mukai and others. We prove that they are smooth 
but very special linear sections of Grassmann varieties, and rational Fano manifolds
of dimension $3a$ and index $a+1$, for $a=1, 2, 4, 8$. We study their 
maximal linear spaces and prove that through the general point pass exactly 
three of them, a result we relate to Cartan's triality principle. We also prove
that they are compactifications of affine spaces.  
\end{abstract}

\maketitle


\section{Introduction}

\subsection{Preliminary:  Reductions of Quadrics}
Let ${\PP}^n = {\PP}(V_{n+1})$, $V_{n+1} = {\CC}^{n+1}$ 
be the complex projective $n$-space, and let 
$\hat{\PP}^n = {\bf  P}(\hat{V}_{n+1})$, 
$\hat{V}_{n+1} = Hom(V_{n+1},{\CC})$, 
the dual projective $n$-space.
Let $Q$ be a smooth quadric in ${\PP}^n$. 
A  {\it non-singular reduction} of $Q$ is any $n$-simplex 
$\Delta \subset {\PP}^n$ (equivalently, any set of $n+1$ independent point in $\hat{\PP}^n$),
 such that in homogeneous 
coordinates $(x) = (x_1:\cdots :x_{n+1})$, defining the $(n-1)$-faces 
$\Delta_i = (x_i = 0)$, $i = 1, ... , n+1$ of $\Delta$, 
the quadratic form $Q(x)$ defining the quadric 
$Q = (Q(x) = 0)$ becomes diagonal. 
The hyperplanes $h_i = (x_i = 0)$, 
spanned on the $(n-1)$-faces $\Delta_i$,
$i = 1,...,n+1$ of $\Delta$ , are the principal axes of 
$Q$ defined by the reduction $\Delta$.

The $(n+1)$ principal axes $h_i$ 
defined by a (nonsingular) reduction $\Delta$ of $Q$ are also 
the vertices of the  dual simplex  $\hat{\Delta} \subset \hat{\PP}^n$ 
to $\Delta$;  and the fact that $\Delta$ is a reduction of $Q$ 
means that the point $Q \in {\PP}(Sym^2 \ \hat{V}_{n+1})$ 
lies in the projective $n$-space ${\Pi}^n = {\Pi}^n({\Delta})$ 
spanned on the Veronese images $\varepsilon_i$ of the axes 
$h_i$, $i = 1,...,n+1$. 
That is: 

\medskip

${\bf ( \ast )}$
{\it The family $Red^o(Q)$ of non-singular reductions of 
the smooth quadric $Q \subset {\PP}^n$ 
is isomorphic to the family of $(n+1)$-secant 
$n$-spaces to the Veronese $n$-fold  $v_2(\hat{\PP}^n)$
passing thourgh $Q$.}

\subsection{Reductions in Jordan algebras.}
We shall see that the projectivized simple  Jordan algebras 
are the natural projective representation spaces where one 
can define analogs of reductions in a way similar to the 
case of quadrics. 

On the one hand, 
the observation is that the space ${\PP}(W) = {\PP}(Sym^2 \
\hat{V}_{n+1})$   is an irreducible projective representation space of the 
group $G = SL_{n+1}$.  Moreover  ${\PP}(Sym^2 \ \hat{V}_{n+1})$ is a
prehomogeneous  projective space of $SL_{n+1}$, 
i.e. the group $G = SL_{n+1}$ acts transitively over an open subset 
of this space  --  the set 
${\PP}(Sym^2 \ \hat{V}_{n+1}) - Det$ 
of quadrics of rank $n+1$. The last identifies the varieties of 
reductions of any two quadrics  of rank $n+1$. 

On the other hand,  
the Veronese variety $v_2(\hat{\PP}^n)$, 
which is the closed orbit of the projective action $\rho$ 
of $SL_{n+1}$ on ${\PP}(Sym^2 \ \hat{V}_{n+1})$, 
is isomorphic to the projective
$n$-space.  This  makes it possible to 
define the reductions as simplices 
with faces defined by intersections of 
coordinate hyperplanes in ${\PP}^n$.

In the above context, it is natural to look 
for smooth projective varieties that are, in some way,  analogs
of the projective space. 
The most natural possible analogs of 
the projective space are the  Severi varieties 
$\AP2$ -- the Veronese models of the projective planes 
over the 4 complexified composition algebras 
$\AA = \RR, \CC, \HH, \OO$,  
and the varieties of Scorza -- 
the projective $n$-spaces $\AA\PP^n$, $n \ge 3$ 
over $\AA = \RR, \CC, \HH$ \cite{roberts,zak}. 

The ambient projective spaces 
of all the varieties $\AA\PP^n$ -- of Severi, as well of Scorza -- 
are projectivized prehomogeneous spaces.
The representation spaces supporting the varieties  
$\AA\PP^n$, $n = 2, a = 1,2,4,8$ and  $n \ge 3, a = 1,2,4$  
are exactly the spaces of all the simple complex Jordan algebras 
${\PP}({\mathcal J}_{n+1}(\AA))$ of rank $n \ge 3$. 
In particular, for $\AA = \RR = {\bf C}$ one obtains again the 
Jordan algebra 
${\mathcal J}_{n+1}(\RR) \cong Sym^2 \ {\bf C}^{n+1}$ of symmetric
matrices of order $n+1$. 

\subsection{The principal results and structure of the paper}  

In this paper we treat the case $n = 3$, i.e. the question about
the  description of the varieties of reductions in the four projectivized 
simple complex  Jordan algebras of order $3$.

The four complex composition algebras $\AA = \RR, \CC, \HH, \OO$ 
are the complexifications of the four normed division algebras 
${\bf A = R,C,H,O}$ -- the reals, the complexes, the quaternions, 
and the octonions. As complex vector spaces these algebras 
have correspondingly dimensions $a = \dim_{{\bf C}} \AA = 1,2,4,8$; 
in particular $\RR = {\bf C}$.

For any $a = 1,2,4,8$, the Jordan algebra $\JA$ is the complex 
vector space of $\AA$-Hermitian matrices of order $3$. Its projectivization 
contains three types of matrices, depending on the {\it rank} (which can be defined
properly even over the octonions). In particular, the (projectivization of  
the) set of rank one matrices is the Severi variety $X_a$, the projective 
$\AA$-plane, a homogeneous variety of dimension $2a$. 

For a point $w\in \PJA$ defined by a rank three matrix, 
a non-singular reduction of $w$ is a $3$-secant plane 
to $X_a$ through $w$ -- see ${\bf 1.1 ( \ast )}$.
The projection from the fixed point $w$ sends the 
quasiprojective set $Y_a^o$ of non-singular reductions of $w$ 
isomorphically to the family of $3$-secant lines 
to the projected Severi variety 
$\overline{X}_a$ inside the projective space $\PJA_o  = {\PP}^{3a+1}$ of
traceless matrices;
and the projective closure $Y_a$ of $Y_a^o$ 
in the Grassmannian $G(2,\JA_o) = G(2,3a+2)$ is {\it the 
variety of reductions} of $w$, our main object of study. 

\smallskip
In the first section we relate the family $Y_a^0$ of simple reduction 
planes with another series of homogeneous varieties $X^a$. These varieties 
appear in the $\CC$-column of the geometric Freudenthal square explored
in \cite{LMfreud} (while the Severi varieties are those of the $\CC$-line). 
Specifically, we note that the choice of $w$ gives an embedding in these
varieties $X^a$ of a copy $X^a_w$ of a variety from the $\RR$-column of
the magic square. Then, we prove that the choice of a reduction plane
gives an embedding in $X^a_w$ of what we call a {\it triality variety} $Z_{\e}$, 
a variety from the $\underline{0}$-column of the magic square. We give an
interpretation of these varieties as zero-set of sections of homogeneous 
vector bundles whose spaces of global sections are precisely the Jordan 
algebras $\JA$. This is another example of the fascinating geometry related
to Freudenthal's magic square, the new feature here being, while we usually 
understand the square line by line, the geometry of the first three columns
is deeply interwoven with that of the first two lines.

\smallskip
In the second section, we focus on the beautiful geometry of the completion $Y_a$ of 
$Y_a^0$, the variety of reductions. This subvariety of $G(2,\JA_0)$ has dimension 
$3a$, and is endowed with a 
natural action of the automorphism group $\SOA:={\rm Aut}(\JA)$ of the Jordan algebra. 
We prove that $Y_a$ has four $\SOA$-orbits, which we describe explicitely: they have
codimension $0$, $1$, $2$ and $4$ (Proposition \ref{orbits}). We prove that $Y_a$ can
be defined as a linear section of the ambient Grassmannian $G(2,\JA_0)$, in its 
Pl\"ucker embedding (Proposition \ref{nontransverse}). For $a>1$, this section is non 
transverse, not even proper. Nevertheless, we prove that $Y_a$ is a smooth subvariety
of $G(2,\JA_0)$ (Theorem \ref{smoothness}). This unexpected phenomenon is related to the 
presence in $Y_a$ of large linear spaces: namely, $Y_a$ is covered by a family of 
$\PP^a$'s parametrized by the Severi variety $X_a$. We prove that the stabilizer 
of a generic point of $Y_a$ is the semi-direct product of a triality group by a 
symmetric group ${\mathfrak S}_4$ (Propositions \ref{s3}). Moreover, exactly three 
$\PP^a$'s pass through that point, which are permuted by the symmetric 
group ${\mathfrak S}_3$, obtained as the quotient of ${\mathfrak S}_4$
by the normal subgroup of permutations given by  products of two disjoint
transpositions. This 
leads to a very nice geometric picture of the Lie algebra isomorphism 
$$\fsoa=\ft(\AA)
\op\AA_1\op\AA_2\op\AA_3$$ (see {\bf 2.7.1} for the notations), which was used in \cite{LMdel}
in a very different context. This geometric occurence of triality completes the 
picture given by E. Cartan in his paper on isoparametric families of hypersurfaces, 
the first geometric appearance of the exceptional group $F_4$ \cite{cartan}. 

\smallskip
Using the geometry of linear subspaces on $Y_a$, we prove that the point-line incidence
variety $Z_a$ over $Y_a$ is the blow-up of the projected Severi variety $\overline{X_a}$
in $\PP\JA_0$ (Proposition \ref{blowup}). An easy consequence is that $Y_a$ is a smooth
Fano manifold of index $a+1$ (and dimension $3a$) with a cyclic Picard group. We also
compute its Betti numbers, and its degree with respect to the Pl\"ucker embedding. 
Finally, we use the group action to prove that $Y_a$ is a minimal compactification 
of $\CC^{3a}$. More precisely, the maximally degenerate hyperplane sections of $Y_a$
are parametrized by its closed orbit (which identifies with the space of special lines
on the hyperplane section $X_a^0$ of the Severi variety), and their complements in 
$Y_a$ are affine cells (Theorem \ref{affine}). Remember that the only minimal projective 
compactification of $\CC^2$ is the projective plane, and that there exists only four types
of such compactifications of $\CC^3$. Several people asked for the classification of 
minimal compactifications of $\CC^n$, but few explicit examples seem to be known. 
Our varieties of reductions give a series of such examples.

\smallskip
The case $a=1$ is classical: $Y_1$ is a transverse intersection of the Grassmannian
$G(2,5)\subset\PP^9$ with a codimension three linear subspace. This is the Fano threefold
of degree $5$, studied in particular by Mukai \cite{mukai}. The fact that $Y_1$ is a 
compactification of $\CC^3$ was discovered by Furushima (see \cite{furu} and the 
references therein). We realized that the second variety of reductions $Y_2$ also
appears in the litterature, in a slightly disguised form: it is called in \cite{tjotta}
the variety of determinantal nets of quadrics. In our interpretation, it is rather
the space of abelian planes in $\fsl_3$, and is therefore closely connected with the 
commuting variety of $\fsl_3$. We show in the last section 
that it is the image of one of the two extremal 
contractions of the punctual Hilbert scheme ${\rm Hilb^3}\PP^2$ (the other one is the 
Hilbert-Chow morphism, whose image is of course singular). Since $Y_2$ is 
a Fano manifold of index $3$ and dimension $6$, a generic codimension three linear
section is a smooth Calabi-Yau threefold. We conclude the paper by a computation 
of the Euler number of this Calabi-Yau manifold, showing that its general deformation 
is not induced by a deformation of the section. 



\bigskip

\section{Reductions in the Jordan algebras $\JA$}

\subsection{The four special complex Jordan algebras of order 3}
${   }$
\hfill \hfill \hfill 
\linebreak
Let ${\bf A = R,} {\bf C,H,O}$, the reals, the complexes, the quaternions 
and the octonions, be the four real division algebras.  
As real vector spaces, they have dimensions correspondingly 
$a = 1,2,4,8$. They are endowed with non degenerate quadratic forms $q$ such that
$q(xy)=q(x)q(y)$. We denote by $q(x,y)$ the associated scalar product, so that 
$q(x)=q(x,x)$. By orthogonal symmetry with respect to the unit element $1$, a vector 
$x \in {\bf A}$ is transformed into its  conjugate $\bar{x} \in {\bf A}$, which is 
such that $x\bar{x}=\bar{x}x=q(x)1$. The real and imaginary parts of $x$ are defined by 
the identities $x=Re(x)1+Im(x)$, $\bar{x}=Re(x)1-Im(x)$.

For any $a = 1,2,4,8$, let 
let $\AA = {\bf A} \otimes_{\bf R} {\bf C}$ 
be the complex composition algebra with multiplication 
$(x_1 \otimes c_1)(x_2 \otimes c_2) = x_1x_2 \otimes c_1c_2$
and conjugation $\overline{x \otimes c} = \bar{x} \otimes c$.
In particular $\RR \cong {\bf C}$, ${\CC} \cong {\bf C} \oplus {\bf C}$, 
and $\HH \cong M_2({\bf C})$, the algebra of 
complex matrices of order $2$. 

For any $a = 1,2,4,8$ the space 
$$\JA = \Bigg\{ 
\begin{pmatrix} 
c_1            & x_3            & \bar{x}_2  \\ 
\bar{x}_3 & c_2            &  x_1            \\ 
x_2            & \bar{x}_1 &  c_3            \\ 
\end{pmatrix}
: c_i \in {\bf C}, x_i \in \AA \Bigg\} \cong {\bf C}^{3a+3}
$$
of $\AA$-Hermitian matrices of order $3$, together with 
the Jordan multiplication 
$A \circ B = \frac{1}{2}(AB+BA)$ 
is a Jordan algebra, i.e. $(\JA,\circ)$ is commutative 
and the equality 
$$(A\circ B) \circ (A \circ A)   = A \circ (B \circ (A \circ A))$$
holds for any $A,B \in \JA$.

\subsection{The Severi varieties and reductions}

On $\JA$ there is a well defined determinant $\det$, a cubic form 
which can be defined in terms of the trace of a matrix and its second
and third powers, by the formula which is usual in $M_3({\bf C})$. 

For $a = 1,2,4,8$ the subgroup  of $GL_{\bf C}({\JA})$ of complex-linear transformations 
preserving  the determinant, is the product of its center by the derived group $SL_3(\AA )$. 
This semi-simple group is isomorphic correspondingly to 
$SL_3, SL_3 \times SL_3, SL_6, E_6$;  and in fact  the action
$\rho_a: SL_3(\AA) \times \JA \rightarrow \JA$
is an irreducible representation of 
$SL_3(\AA)$ in the complex vector space 
$\JA = {\bf C}^{3a+3}$. 
More precisely, let $A_2, A_2 + A_2, A_5, E_6$ 
be correspondingly the Dynkin diagrams of 
$SL_3(\AA)$, $a = 1,2,4,8$. 
Then, in the notation of [Bourbaki], 
we may consider that: 

\medskip

-- The representation $\rho_1$ 
is defined by the weight $2\omega_1$, 
i.e. $\rho_1$ is the 2-nd symmetric power of the 
standard representation of $SL_3$; 
in particular the Jordan algebra $\JR$ is isomorphic  
to the algebra $Sym^2 \ {\bf C}^3$ of symmetric 
complex matrices of order $3$.   

-- The representation $\rho_2$ 
is defined by any pair of weights 
$(\om_i',\om_j'')$, $1\le i,j\le 2$
of the two copies $(A_2',A_2'')$ 
of $A_2$ in the diagram of $SL_3 \times SL_3$; 
in particular $\JC$ is isomorphic to the 
algebra $\otimes^2 \ {\bf C}^3$ of 
matrices of order $3$. 

-- The representation 
$\rho_4$ is defined by $\omega_2$ or by $\omega_4$ , 
i.e. $\rho_4$ is the second or the fourth alternative power 
of the standard representation of $SL_6$; 
in particular $\JC$ is isomorphic to the 
algebra $\wedge^2 \ {\bf C}^6$ of antisymmetric 
matrices of order $6$. 

-- The representation $\rho_ 8$ 
is defined by any of the weights $\omega_1$ or 
$\omega_6$  of $E_6$. 

\medskip 

Denote by $\rho_a$ also the projectivized action 
$$\rho_a: SL_3(\AA) \times \PP\JA \rightarrow \PP\JA$$
on the projective complex space $\PP\JA = {\PP}^{3a+2}$.
The following is well-known, see e.g. \cite{sk,LMfreud}:

\begin{lemm}\label{sevorbits}
The projective action of $\rho_a$ of $SL_3(\AA)$
splits  $\PP\JA$ into a union of $3$ orbits 
$$
\PP\JA = (\PP\JA - \DA) \ \cup   \ (\DA - \AP2) \ \cup \ \AP2,$$
where $\DA$ = (det = 0) and 
$\AP2$ are correspondingly the determinantal cubic hypersurface 
and the locus of Jordan $\AA$-matrices of order $3$ and of rank $1$. 
Moreover $\DA = Sec(\AP2)$ = the union of all the secant 
lines to $\AP2$; and $\AP2 = Sing  \ \DA$. 
\end{lemm}

\medskip
The {\it rank} of a matrix with coefficients in $\OO$ is a rather delicate notion. 
In this paper, it will suffice to consider that by definition, these three orbits
in $\PP\JA$ consist in matrices of rank $3$, $2$ and $1$ respectively. 

\smallskip
The varieties $\AP2 = \RP2, \CP2, \HP2$ and $\OP2$ 
can be interpreted as the four complex projective 
$\AA$-planes, see \cite{roberts}. From another point of view 
they are also all the four Severi varieties \cite{roberts, zak, chaput}. 
They can be described very explicitely:

\begin{lemm} The Severi variety $X_a\subset\PP\JA$ is defined by the
equations $X^2=\trace(X)X$, which generate its ideal. Its intersection 
with the affine subspace of $\PP\JA$ on which the first diagonal coefficient is non 
zero is
$$X_a\cap\{c_1\ne 0\}=\Bigg\{ 
\begin{pmatrix} 
1           & x            & y  \\ 
\bar{x} & \bar{x}x         &  \bar{x}y \\ 
\bar{y}  & \bar{y}x & \bar{y}y
\end{pmatrix}, \quad  x,y,\in\AA\Bigg\}\cong {\bf C}^{2a}.$$
\end{lemm}

The Severi varieties really show up the geometry of projective planes. 
They are covered by a family of projective $\AA$-lines $\AA\PP^1\cong \QQ^a$, 
quadrics of dimension $a$, and two (generic) such $\AA$-lines intersect at a unique 
point. 

\smallskip   
In the interpretation $\JR = Sym^2 \ {\bf C}$, the 
$\RR$-plane $\RP2$ is the Veronese image $v_2({\PP}^2)$
of the complex projective plane, and the hypersurface 
$\DR$ is the symmetric determinant cubic 
in $\PP( Sym^2 \ {\bf C})$. 
For $\JR \cong \otimes^2 \bf C^3$, the $\CC$-plane 
$\CP2$ is the Segre variety ${\PP}^2 \times {\PP}^2$, 
and  $\DC$ is the symmetric determinant cubic 
in $\PP( \otimes^2 \ {\bf C})$. 
For $\JH \cong \wedge^2 \ {\bf C}^6$, the complex 
quaternionic plane $\HP2$ is the grassmannian $G(2,6)$, 
and $\DH$ is the pfaffian cubic 
$Pf \subset \PP(\wedge^2 \ {\bf C}^6)$. 
At the end, the complex octonionic plane, 
or the complex Cayley plane $\OP2 \subset \PP(\JO ) ={\PP}^{26}$ 
is a smooth Fano $16$-fold of degree $78$
with $Pic \ \OP2 = \ZZ H$, \ $H$ being the hyperplane section, 
and $K_{\OP2} = -12H$.   
\begin{center}
\setlength{\unitlength}{3mm}
\begin{picture}(40,4)(-2,-2)

\put(0,0){$\circ$}
\put(2,0){$\bullet$}
\put(0.5,.4){\line(1,0){1.6}} 
\put(2,1){{\small $2$}}

\multiput(6,0)(4,0){2}{$\bullet$}
\multiput(8,0)(4,0){2}{$\circ$}
\multiput(6.5,.4)(4,0){2}{\line(1,0){1.6}} 
\put(8.7,0){$\times$}

\multiput(16,0)(2,0){5}{$\circ$}
\multiput(16.5,.4)(2,0){4}{\line(1,0){1.6}} 
\put(18,0){$\bullet$}

\put(28,0){$\bullet$}
\multiput(28,0)(2,0){5}{$\circ$}
\multiput(28.5,.4)(2,0){4}{\line(1,0){1.6}} 
\put(32.3,-1.45){\line(0,1){1.5}}
\put(32,-2){$\circ$}
\end{picture} 
\end{center}

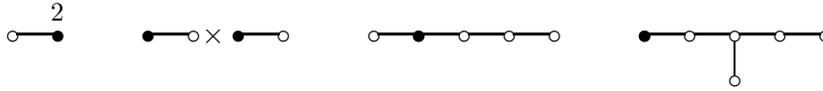
\centerline{{\small\it Weighted Dynkin diagrams of the Severi varieties $X_a$}}

\bigskip
To the four Jordan algebras one can attach the 
algebra $\J0 \cong {\bf C}^3 $ of complex diagonal matrices 
of order $3$, coming from the algebra 
$\underline{0} = (0)$. 
Then the determinant hypersurface 
$\DD(\underline{0}) \subset \J0 
= {\PP}^2$ is, of course, the coordinate triangle 
$\Delta \subset {\PP}^2$, 
and the $\underline{0}$-plane (or the $0$-th Severi variety) 
$\underline{0}{\PP}^2 = \vee^3 \ {\PP}^0$ 
is the triple of vertices of $\Delta$.

The five Severi varieties fill in the $2$-nd line of the 
extended Freudenthal square, see [LM]:  

\
 
\begin{center}\begin{tabular}{|c|c|cccc|c} \cline{1-6} 
 & \underline{0} & 
 ${\RR}$ & ${\CC}$  & ${\HH}$ & ${\OO}$ \\ \cline{1-6} 
 
${\RR}$ & $\emptyset$ & 
$v_4({\PP}^1)$ & ${\PP}(T_{{\PP}^2})$ & $IG(2,6)$ & 
${\OO}{\PP}^2_{0}$ &  
{\it section of Severi} \\ 
 
${\CC}$ & $\vee^3 {\PP}^o$ & $v_2({\PP}^2)$ & 
$\times^2 {\PP}^2$ & $G(2,6)$ & ${\OO}{\PP}^2$ &  
{\it Severi} \\ 
 
${\HH}$ & $\times^3 {\PP}^1$ & 
$LG(3,6)$ & $G(3,6)$ & ...  & ... &  
{\it Legendre} \\ 
 
${\OO}$ & $D_4^{ad}$ & 
$F_4^{ad}$ & $E_6^{ad}$ & ...  & ... &  
{\it adjoint} \\ \cline{1-6} 

\end{tabular}\end{center} 

\

Here we denoted by $IG(2,6)$ (resp. $LG(3,6)$) the isotropic (resp. lagrangian) grassmannian
of isotropic $2$-planes (resp. $3$-planes) in $\CC^6$ with respect to a symplectic form.  
The adjoint varieties are the closed orbits in the projectivizations of the adjoint representations
of the simple Lie algebras. 

\smallskip
As we will see below, the $\underline{0}$-th column 
will play an important role in the description of the 
varieties of reductions in the Jordan algebras 
$\JA, a = 1,2,4,8$. 

\medskip

\begin{defi} {\bf (Reductions and reduction planes)}
Let $a = 1,2,4,8$, and let  $w \in \PJA^o = \PJA - \DA$
be a projective rank $3$ Jordan matrix of order $3$. 
Call a non-singular reduction of $w$ any simply 
$3$-secant plane 
${\PP}^2 \subset {\PP}^{3a+2} = \PJA$
to $\AP2$ which passes through the point $w$.\end{defi}

\subsection{The four varieties of reductions $Y_a \subset G(2,\JA_0)$}

Let $w \in \PJA^o$ = $\PJA - \DA$ be as above, 
and let 
$$
\pi_w:  \PJA \dashrightarrow \PJA_w
$$
be the rational projection from $w$. 
One can identify the base space 
$\PJA_w$ of $\pi_w$ with the 
polar hyperplane 
${\PP}^{3a+1}_w \subset {\PP}^{3a+2}$
to the rank $3$ point $w$,  defined has follows: let $\det(X,Y,Z)$ be the 
polarisation of the determinant, i.e. the unique symmetric trilinear form on 
$\JA$ such that $\det(X)=\det(X,X,X)$. Then the polar hyperplane to $w$ is defined
by the equation $\det(W,W,X)=0$, where $W\in\JA$ is any representative of $w$. 

Since the point $w$ doesn't lie in the secant variety 
$Sec(\AP2) = \DA$ (see Lemma \ref{sevorbits}), then 
the projection ${\pi}_w$ sends:
{\bf (a)} the Severi variety $X_a = \AP2$ isomorphically to 
its image $\overline{X}_a \subset {\PP}^{3a+1}_w$; 
{\bf (b)} any non-singular reduction ${\PP}^2$ 
of $w$ to a line $l$
which is simply $3$-secant to $\bar{X}_a$. 

Inversely, any line $l \subset {\PP}^{3a+1}_w$
which is simply $3$-secant to $\overline{X_a}$
is a projection from $w$ of a unique plane 
${\PP}^2 \in Y^o_w$. 
Therefore $\pi_w$ embeds  the set of reductions 
of $w$ as a subset  $Y^o_{a,w}$ of the grassmannian 
$G(2,\JA_w)$ of lines in ${\PP}^{3a+1}_w=\PP\JA_w$, 
which yields the following

\begin{defi}
For the point $w \in \PJA^o = \PJA - \DA$, 
define the variety $Y_{a,w}$ to be the closure 
of $Y_{a,w}^o$ in the grassmannian 
$G(2,\JA_w)=G(2,3a+2)$. 
\end{defi}
 
Since the group $SL_3(\AA)$ acts transitively on the points 
$w \in \PJA^o$, as well as on the Severi variety $\AP2$ 
then all the varieties $Y_{a,w}$, 
$w \in \PJA^o$ are projectively equivalent, 
by the induced action of $SL_3(\AA)$ on 
$G(2,\JA)$, to the same variety $Y_a$; 
and we let 
$$Y_a := Y_{a,I}$$
where $I$ is the projectivized unit matrix 
$diag(1,1,1) \in \JA$. Note that $\JA_I$ coincides with $\JA_0$, the space of traceless
matrices in $\JA$. 

\subsection{The $\CC$-column}

Denote by 
$X^a$ = $\PP(T_{\PP}^2)$, $\times^2{\PP}^2, LG(3,6)$, $E_6^{ad}$, 
$a = 1,2,4,8$ the four varieties of the ${\CC}$-column of the Freudenthal square.
By the preceding, for any $a = 1,2,4,8$ the choice of the 
Jordan algebra $\JA$ defines uniquely the group 
$\SLA$, and the variety 
$X_a = \AP2 \subset \PJA$ in  the $\CC$-line 
of the Freudenthal square. 
Similarly, $\JA$ defines uniquely the variety 
$X^a$ in the $\CC$-column 
as the closed orbit of the projective representation 
$\rho^a$ of the group $\SLA$, as follows; 
as above we use the notation of weights from \cite{bou}: 
 
\medskip

-- $\rho^1$ is the adjoint representation of $SL_3$, 
defined by the weight $\om_1 + \om_2$
of $A_2$.

-- $\rho^2$ 
is defined by any pair of weights 
$(\om_i',\om_j'')$, $1\le i,j\le 2$
of the two copies $(A_2',A_2'')$ 
of $A_2$ in the diagram of $SL_3 \times SL_3$; 
in particular $\JC$ is isomorphic to the 
algebra $\otimes^2 \ {\bf C}^3$ of 
matrices of order $3$. 

-- $\rho^4$ is defined by the weight $\omega_3$ 
of $A_5$ 
giving the $3$-rd alternative power 
of the standard representation of $SL_6$; 

-- $\rho^8$ is the 
adjoint representation of the group $E_6$ 
given by the weight $\omega_2$ of the diagram $E_6$.

\subsection{Points in $\PJA$ and isotropic varieties from the
$\RR$-column}

\subsubsection{Points and isotropic groups}

Let $\rho_a$ be the action $\SLA$ in $\PJA$. For a point $w \in
\PJA^o$  define $\SOA_w \subset \SLA$ to be the connected component in the 
isotropy subgroup of $w$. 
Since $\PJA^o$ is a $\rho_a$-orbit of $\SLA$, then  the group $\SLA$ 
permutes the set 
$$\{ \SOA_w \subset \SLA , w \in \PJA^o \} $$
of all these copies this way identifying any of them with the 
subgroup $\SOA = \SOA_I \subset \SLA$ preserving the projective unit 
matrix $I = diag(1,1,1)$. By \cite{chaput2}, Proposition 3.2, this group
$\SOA$ coincides with the automorphism group ${\rm Aut}\JA$ of the 
Jordan algebra $\JA$. It preserves not only the determinant, but also the 
linear form $\trace(X)=\det(I,I,X)$ and the quadratic form $Q(X)=\trace(X^2)
=\det(I,I,X)^2-2\det(I,X,X)$. 

The following result is well-known, see e.g. \cite{ahiezer}:

\begin{lemm}\label{2orbits}
 The action of $\SOA$ on the Severi variety $X_a$ has exactly 
two orbits, the hyperplane section $X_a^0=X_a\cap\PP\JA_0$ and its complement.
\end{lemm}

For $a = 1,2,4,8$,  $\SOA$ is correspondingly
$SO_3$, $SL_3$, $Sp_6$, $F_4$,   see e.g. \cite{LMfreud}. We 
denote by $SO_{3,w}$, $SL_{3,w}$, $Sp_{6,w}$, $F_{4,w}$ 
the copy of $\SOA_w$ in 
$\SLA_w$ = $SL_3$, $SL_3 \times SL_3$, $SL_6$, $E_6$ 
defined by $w \in \PJA^o$. 

\medskip
For a fixed, $w \in \PJA^o$, the projective representation $\rho^a$ 
of $\SLA$ induces a projective representation $\rho^a|_{w}$
of the subgroup $\SOA$. The representation $\rho^a|_{w}$
is already reducible; in fact one has, in terms of highest weights:
\medskip

-- for $a = 1$,  \ $\rho^1|_{w} = 4\om_1 \oplus 2\om_1$, 
and we let $\rho^{1}_{w} = 4\om_1$.  

-- for $a = 2$, \ $\rho^2|_{w} = 2\om_1 \oplus \om_1$, 
and we let $\rho^{2}_{w} = 2\om_1$.  

-- for $a = 4$, \ $\rho^4|_{w} = \om_3 \oplus \om_1$
and we let $\rho^{4}_{w} = \om_3$.

-- for $a = 8$, \ $\rho^8|_{w} = \om_1 \oplus \om_4$,
and we let $\rho^{8}_{w} = \om_1$. 
\medskip

For any $a = 1,2,4,8$, the subrepresentation 
$\rho^{a}_{w}$ of $\SOA_w$ is irreducible; 
and the choice of $w \in \PJA$ defines 
uniquely the projective representation subspace
$\PP(V^a_w)$ of $\rho^{a}_{w}$, in the space $\PP(V^a)$ 
of $\rho^a$. For $a = 1,2,4,8$  denote these subspaces 
correspondingly by 
$$\begin{array}{cc}
\PP(Sym^4_w{\bf C}^2) \subset \PP(sl_3), &
\PP(Sym^2_w{\bf C}^3) \subset \PP(\otimes^2 \ {\bf C}^3), \\
\PP(\wedge^{(3)} _w{\bf C}^6) \subset \PP(\wedge^{3} _w{\bf C}^6), &
\PP(\ff_{4,w}) \subset \PP(\fe_6). 
\end{array}$$

In particular $dim \ \PP(V^a_w) = 4, 5, 13, 51$ 
in the space $\PP(V^a)$ of dimension  $7, 8, 19, 77$.

\begin{defi}
{\bf (Isotropic spaces and isotropic varieties)}     
For $w \in \PJA^o$, $a = 1,2,4,8$ we call 
the subspace $V^a_w \subset V^a$ 
(respectively the projective subspace 
$\PP(V^a_w) \subset \PP(V^a)$) 
the isotropic subspace of $w$
(respectively the isotropic projective 
space of $w$).  
We call the closed $\SOA_w$-orbit 
$X^a_w \subset \PP(V^a_w)$
the isotropic variety of $w$.  
\end{defi}

For $a = 1,2,4,8$, denote the isotropic subvariety 
$X^a_w \subset X^a$ by 
$$
v_4({\PP}^1)_w \subset \PP(T_{\PP}^2), \quad 
v_2({\PP}^2)_w \subset \times^2 \PP^2, \quad
LG_w(3,6) \subset G(3,6) \quad \mbox{ and } 
F_{4,w}^{ad} \subset E_6^{ad}.$$ 
Now the following is direct:

\begin{lemm}\label{isotropic}
Let $a = 1,2,4,8$. Then any $w \in \PJA^o$ defines uniquely 
the isotropic subvariety $X^a_w \subset X^a$ in the 
$\RR$-column of the Freudenthal square as the intersection 
$$
X^a_w = X^a \cap \PP(V^a_w)
$$
of  $X^a \subset \PP(V^a)$ with the isotropy projective subspace 
$\PP(V^a_w)$ of $w$. 
\end{lemm}

\begin{center}
\setlength{\unitlength}{3mm}
\begin{picture}(40,13)(-3,-2)

\put(1,4){$\bullet$}
\put(1,3){{\small 4}}
\put(0,0){$\bullet$}
\put(2,0){$\bullet$}
\put(0.5,.4){\line(1,0){1.6}} 

\put(8.7,4){$\bullet$}
\put(8.95,4.6){\line(0,1){1.6}}
\put(8.7,6){$\circ$}
\put(8.7,3){{\small 2}}
\multiput(6,0)(2,0){4}{$\circ$}
\multiput(8,0)(2,0){2}{$\bullet$}
\multiput(6.5,.4)(4,0){2}{\line(1,0){1.6}} 
\put(8.7,0){$\times$}

\multiput(16,0)(2,0){5}{$\circ$}
\multiput(16.5,.4)(2,0){4}{\line(1,0){1.6}} 
\put(20,0){$\bullet$}
\multiput(20,6)(0,2){2}{$\circ$}
\put(20,4){$\bullet$}
\put(20.25,6.6){\line(0,1){1.5}}
\multiput(20.15,4.4)(.3,0){2}{\line(0,1){1.7}}

\multiput(28,0)(2,0){5}{$\circ$}
\multiput(28.5,.4)(2,0){4}{\line(1,0){1.6}} 
\put(32.3,-1.45){\line(0,1){1.5}}
\put(32,-2){$\bullet$}
\multiput(32,6)(0,2){3}{$\circ$}
\put(32,4){$\bullet$}
\put(32.25,4.6){\line(0,1){1.6}}
\put(32.25,8.5){\line(0,1){1.6}}
\multiput(32.15,6.5)(.3,0){2}{\line(0,1){1.55}}
\end{picture} 
\end{center}

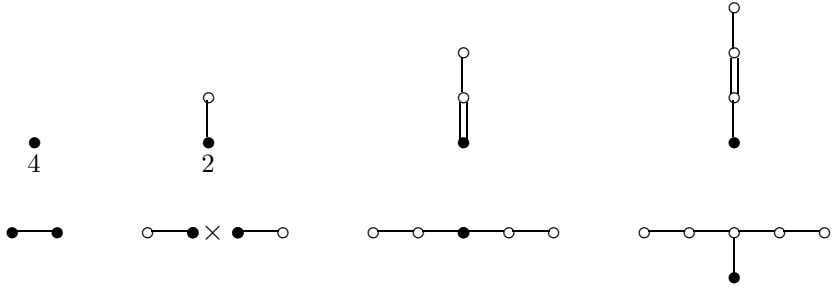
\centerline{{\small\it Weighted Dynkin diagrams of the varieties $X^a$ and $X^a_w$}}

\medskip
The weighted Dynkin diagrams of the varieties $X_a^w$ are obtained by {\it folding}
those of the varieties $X_a$. 

\subsubsection{\bf Isotropic varieties as zero-sets}

In this section we let $a>1$.
It turns out that the isotropy subvarieties 
$X^a_w \subset X^a$ are zero-sets of sections of a homogeneous 
vector bundle ${\cE}^a$ on $X^a$. 
Recall that an irreducible homogeneous vector bundle on a homogeneous variety
$X=G/P$ is determined by the highest weight of the corresponding $P$-module,
which can be encoded in a weighted Dynkin diagram. To get the weighted Dynkin
diagram of our vector bundle on $X^a$, we just superimpose the weighted diagram 
of the Severi variety $X_a$, to that of $X^a$. 

Since the weighted Dynkin diagram of $X^a$ has a twofold symmetry, while that 
of $X_a$ does not, we obtain in fact {\it two} vector bundles 
$\cE^a$ and $\cE_*^a$
on $X^a$, which can be deduced one from the other through the action of an outer 
involutive automorphism of $X^a$. 

\begin{prop}  
The vector bundles $\cE^a$ and $\cE^a_*$ on $X^a$ are generated by their global 
sections, and their spaces of global sections are isomorphic to the Jordan algebra
$\JA$ and its dual $\JA^*$. Their ranks and determinants are as follows:
$$\begin{array}{rcl} 
rank \ \cE^8 = 6, & \quad & det \ \cE^8 = \cO(3), \\
rank \ \cE^4 = 3, & &  det \ \cE^4 = \cO(2), \\
rank \ \cE^2 = 2, & &  det \ \cE^2 = \cO(1,2).
\end{array}$$
\end{prop}

\proof The first assertion is an immediate consequence of the Borel-Weil theorem. 
The rank of $\cE^a$ (and $\cE^a_*$) can be read off its weighted Dynkin diagram, 
since the $P$-module that defines $\cE^a$ is encoded in the weighted diagram which
is  obtained after deleting the black nodes that define $X^a$. To compute the 
determinant, we need to list the weights of this $P$-module, which are just the
images of the highest weight by the Weyl group of $P$ (indeed, theses modules
are minuscule, as we can see case by case). Taking the sum of these weights, 
we get the weight of the determinant. \qed

\medskip
More explicitly, the vector bundles $\cE^a$ and $\cE^a_*$ on $X^a$ can be 
described as follows. On $X^2=\PP^2\times\PP^2$, we have the pull-backs $\cO(1)$ 
and $\cO(1)'$ of the hyperplane line bundles on the two copies of $\PP^2$, and
the pull-backs $T$ and $T'$ of the rank two tautological bundles; then $\cE^2$
and $\cE^2_*$ are the bundles $Hom(T,\cO(1)')$ and $Hom(T',\cO(1))$. 
On $X^4=G(3,6)$, let $T$ and $Q$ denote the tautological and quotient vector bundles, 
both of rank three; then $\cE^4$ and $\cE^4_*$ are the bundles $\wedge^2Q$ and $\wedge^2T^*$. 

The case of $X^8=E_6^{ad}$ is slightly more subtle. Recall that the adjoint variety 
$E_6^{ad}$ is the closed orbit in $\PP\fe_6$. If $X\in\fe_6$ defines a point of the
adjoint variety, consider its action $d\rho_8(X)$ on $\JO$ (see {\bf 2.2}). We claim
that $d\rho_8(X)$ has rank six. To check this, we use the fact that $\fe_6$ contains 
a copy of $\fso_8$ whose action on $\JO$ can be described very explicitely. Recall
that the {\it infinitesimal triality principle} asserts that for any $g=g_1\in\fso_8$, 
there exists uniquely defined operators $g_2,g_3\in\fso_8$ such that $g_2(xy)=xg_1(y)
+g_3(x)y$ for all $x,y\in\OO$. Then the action of $g$ on $\JO$ is given by the formula
\cite{jacobson,harvey,LMpop}
$$d\rho_8(g)
\begin{pmatrix} 
c_1            & x_3            & x_2  \\ 
\bar{x}_3 & c_2            &  x_1            \\ 
\bar{x_2}            & \bar{x}_1 &  c_3            
\end{pmatrix}
=
\begin{pmatrix} 
0            & g_3(x_3)            & g_2(x_2)  \\ 
\overline{g_3(x_3)} & 0            &  g_1(x_1)            \\ 
\overline{g_2(x_2)}    & \overline{g_1(x_1)} &  0             
\end{pmatrix}.$$
When $g$ belongs to the adjoint variety of $\fso_8$ (which is naturally contained in that of 
$\fe_6$), $g_1, g_2$ and $g_3$ have minimal rank, that is 
rank two, so clearly $d\rho_8(X)$ has rank six. Therefore we get two vector bundles
whose fibers at $X\in E_6^{ad}$ are $\JO/Ker\ d\rho_8(X)$ and $(Im\  d\rho_8(X))^*$,
respectively.  
These bundles are homogeneous of rank six, and respectively quotients of the trivial
bundles with fibers $\JO$ and $\JO^*$: they are our bundles $\cE^8$ and $\cE^8_*$.

\medskip
Distinguishing the two bundles $\cE^a$ and $\cE^a_*$ is really a matter of convention. 
Our choice will be such that the space of global sections of the vector bundle $\cE^a$ 
is the Jordan algebra $\JA$, rather than its dual. 

Since for a section $w \in H^0(X^a, \cE^a)^o = \JA^o$, 
and $c \in \CC^*$, the zero-sets $Z(w)$ and $Z(cw)$ 
co\"{\i}ncide, one can regard equivalently the elements  
$w \in \PJA$ as {\it projective sections} of $\cE^a$, 
and their zero-sets $Z(w) \subset X^a$.

\begin{prop}
Let $a>1$. Then for any projective section 
$$w \in \PP(H^0(X^a, \cE^a)^o) = \PP\JA^o,$$  
the zero-set $Z(w) \subset X^a$ coincides with the isotropy 
subvariety $X^a_w \subset X^a$. 
\end{prop} 

\proof We treat the case $a=8$, the other ones are simpler. Since $\PP\JO^0$
is an orbit of $E_6$, we may suppose that $w=I$, the identity of $\JO$. By the 
description we have just given of $\cE^8$, a point $X\in E_6^{ad}$ belongs to 
the zero-set $Z(I)$ if and only if $I\in Ker\ d\rho_8(X)$, which means that $X$ 
belongs to the isotropy Lie algebra of $I$. But recall that the (connected component
of the identity in the)  isotropy group of $I$ is the automorphism group ${\rm Aut}
\JO=F_4$, hence the isotropy Lie algebra is $\ff_4$. We conclude that $Z(I)=
E_6^{ad}\cap\PP\ff_4=F_4^{ad}$. \qed

\subsection{The 3-secant Lemma}

Let $G(2,\JA)$ be the grassmannian of lines in $\PJA$,
and let $\Delta(\AA) \subset G(2,\JA)$ be the subset 
of lines $L \subset \PJA$ which are not simply 3-secant to 
the determinant hypersurface $\DA$. Clearly
$\Delta(\AA)$ is a hypersurface in $G(2,\JA)$, 
so $\Delta(\AA) \in |{\cO}_{G(2,\JA)}(d)|$ for some $d$
in the Pl\"ucker polarization of $G(2,\JA)$; 
and we shall see that $d = 6$. 
Indeed,  $d$ is the number of intersection points of $\Delta(\AA)$ 
with the general line  $\Lambda \subset G(2,\JA)$. Such a general 
line $\Lambda$ is a plane pencil of lines  in a general plane
${\PP}^2 \subset \PJA$ passing through  a fixed general point 
$w \in {\bf P}^2$. The plane ${\bf P}^2$ intersects $\DA$ at a 
smooth cubic $C$; and now the Hurwitz formula implies that through 
$w$ pass exactly $d = 6$ lines $L_i \in \Lambda$, $i = 1,...,6$ 
which  are tangent to $C$.

\begin{lemm}\label{3secant}
Let $L \subset {\PP}\JA$ be any simply 
3-secant line to $\DA$, i.e. $L \in G(2,\JA)^o$.
 
Then in $\PJA$ there exists a unique simply 3-secant plane 
(or a reduction plane -- see {\bf 1.5}) 
to the Severi variety $X_a = \AP2$ which passes through $L$. 
Denote this plane by 
$$
\PP^2 = \PP^2_{\e}  = \langle\e_1, \e_2,\e_3\rangle,
$$ 
$\e_1,\e_2,\e_3$ being the three intersection points of
$\PP^2$ with $X_a$. The intersection 
$$
\Delta_{\varepsilon} = {\PP}^2_{\varepsilon} \cap D_3({\AA})
$$ 
is a triangle with sides $L_k = \langle\varepsilon_i,\varepsilon_j\rangle$, 
$\{ i,j,k \}  =  \{ 1,2,3 \}$. 
\end{lemm} 

\proof
Let $\gamma_1,\gamma_2,\gamma_3$ be the three intersection 
points of $L$ and $\DA$. By \cite{zak}, for any 
$\gamma_i$, $i = 1,2,3$ there exists a unique subspace 
${\PP}^{a+1}_{i} \subset \PJA$ which passes through $\gamma_i$ and intersects 
the Severi variety $X_a$ along a smooth $a$-dimensional 
quadric $Q_i \subset {\PP}^{a+1}_i$. 
Moreover any space $\PP^{a+1}_i$ is swept out by all the 
lines through $\gamma_i$ which are bisecant or tangent 
to $X_a$; and any two quadrics $Q_i$ and $Q_j$  intersect each 
other at a unique point $\e_k$, $\{i,j,k\} =\{1,2,3\}$.

\smallskip
We shall  see that the plane $\langle\e_1,\e_2,\e_3\rangle$ pases through $L$.
Indeed, let e.g. $\e_3\in X_a$ be the intersection point 
of the quadrics $Q_1$ and $Q_2$, let 
$\PP^2_{3}$ be the plane spanned by $L$ and $\e_3$, 
and let $L_i$ be the line $\overline{\gamma_i\e_3}$, $i = 1,2$. 

Since the line $L_i$ lies in $\PP^{a+1}_1$
and passes through $\gamma_i$, it is bisecant 
to $X_a$: the intersection $L_i\cap X_a$ contains $\e_3$ and another point $\e_i'$. 

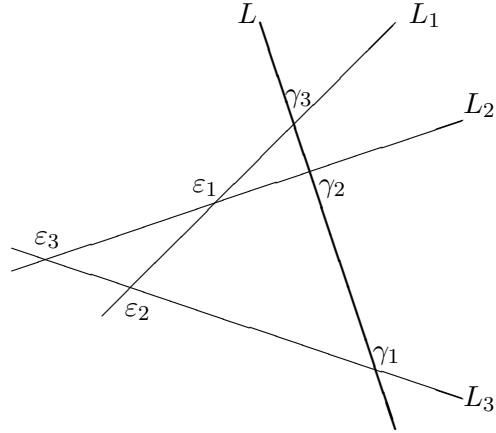
\begin{center}
\setlength{\unitlength}{3mm}
\begin{picture}(40,20)(-6,-6)
\put(0,2){\line(3,1){20}}
\put(0,3){\line(3,-1){20}}
\put(4,0){\line(1,1){13}}
\put(1,3.2){$\varepsilon_3$}
\put(5,0){$\varepsilon_2$}
\put(8,5.5){$\varepsilon_1$}
\put(12.1,9.5){$\gamma_3$}
\put(13.6,5.5){$\gamma_2$}
\put(16,-1.9){$\gamma_1$}
\put(10,13){$L$}
\put(17.6,13){$L_1$}
\put(20,9){$L_2$}
\put(20,-4){$L_3$}

\thicklines
\put(11,13){\line(1,-3){6}}
\end{picture}
\end{center}

Since $\PP^2_{3} \supset L$ and the general point of $L$
has rank $3$, the plane $\PP^2_{3}$ can't lie entirely in 
the determinant cubic  $\DA$. Therefore, it intersects 
$\DA$ along a $1$-cycle $\Delta$ of degree $3$. 
Being bisecant to $X_a$, the  two lines 
$L_1,L_2$ are components of $\Delta$. 
Moreover $\Delta$ passes through the 
rank $2$ point $\gamma_3 \in  \PP^2_3$. 
Therefore, there exists another line $L_3 \subset \PP^2_3$, 
passing through $\gamma_3$, such that the cubic cycle $\Delta = L_1 + L_2 + L_3$. 
This line is bisecant to $X_a$ at its intersection points with $L_1$ and $L_2$, which 
must coincide with $\e_1'$ and $\e_2'$; and now the uniqueness of the 
triple $(\e_1,\e_2,\e_3)$ yields
$\e_i' = \e_i$ and 
$\PP^2_3 = \PP^2_{\e}$.

In addition, if in the above construction one assumes 
that one of the lines, say $L_1$ degenerate to a tangent to 
$X_a$ at $\e_3$, i.e. $\e_2' = \e_3$, then the line $L_2$ 
will degenerate to a multiple (double or triple) component 
of the intersection cubic cycle $\Delta$ of ${\PP}^2_3$,
i.e. either $\Delta = L_1 + 2L_2$, or $L_1 = L_2$ and 
$\Delta = 3L_2$.
But then $L$ can't be simply 3-secant to $\DA$, 
since $L \subset \PP^2_3$ and the intersection cycle 
$L \cap \DA$  must be contained (as a set) in the support of 
$\Delta$. \qed

\begin{coro} \label{trans}
The group $\SLA$ acts transitively  on the set 
of reduction planes in $\PJA$. 
\end{coro}

\proof From the proof of the $3$-secant lemma, we see that
a reduction plane, i.e. a simply $3$-secant plane to $X_a$, is 
uniquely defined by a simply $3$-secant line to $\DA$. So we just need 
to prove that $\SLA$ acts transitively on the set of these lines. 

But it was observed in \cite{sk} that the action of $\SLA$ on $\JA\op\JA$ 
is prehomogeneous, the dense orbit being given by the complement of the 
discriminant hypersurface, that is, precisely,  the set of couples $(x,y)$ 
in $\JA\op\JA$ which are independent, and such that the line $\overline{xy}$
is simply trisecant to $\DA$. Our claim obviously follows. \qed 

\medskip\noindent {\it Remark.} 
Of course, the prehomogeneity of $\JA\op\JA$ is not an accident. Indeed, 
consider the exceptional Lie algebras $\ff_4, \fe_6, \fe_7, \fe_8$. Their 
adjoint representations correspond to an extremal node of their Dynkin diagram. 
Delete its unique neighbour. Then the resulting (non connected) Dynkin diagram
is that of $\SLA\times SL_2$, and if we darken the nodes which neighboured the 
one which has been erased, we get the weighted Dynkin diagram of its representation 
$\JA\otimes\CC^2=\JA\op\JA$. By a Theorem of Vinberg, in such a situation, the action 
of $\SLA\times GL_2$ on $\JA\op\JA$ is prehomogeneous (see \cite{ruben}). 
In a sense, the previous Corollary is thus really a consequence of the existence of 
exceptional simple Lie algebras.

\subsection{Lines in $\PJA$, reduction planes 
                        and triality varieties from the  $\00$-column}

\subsubsection{Lines, reduction planes and triality subgroups} 

In the hypotheses and notations of Lemmas \ref{isotropic} and \ref{3secant}, 
we shall see that the isotropic varieties $X^a_w$
of all the points $w \in \PP^2_{\e} - \Delta$
have a common subvariety -- a {\sl triality 
subvariety} $Z_{\e}$.  We shall regard 
one of the rank 3 points $w \in L$ as fixed; 
this way identifying the unique $3$-secant plane
through $L \supset w$ to a non-singular reduction
of $w$. Let $v \in L\subset \PJA^o$ be another point of $L$
and let $\Delta = \PP^2_{\e} \cap \DA$.

\begin{lemm}
The isotropy subgroup $T^a_{w,v} = \SOA_v \cap \SOA_w$ of the pair 
$(w,v)$ is, up to a finite group, a copy of the triality group $T(\AA)$. 
\end{lemm}

Recall from \cite{LMpop} that this triality group is defined as 
$$T(\AA)=\{g=(g_1,g_2,g_3)\in SO(\AA)^3, \; g_2(xy)=g_1(x)g_3(y)\; \forall x,y\in\AA\}.$$
We have $T(\AA)=1,G_m \times G_m, SL_2 \times SL_2 \times SL_2, Spin_8$ respectively for 
$a=1,2,4,8$. Denote its Lie algebra by 
$$\ft(\AA)=\{u=(u_1,u_2,u_3)\in \fso(\AA)^3, \; u_1(xy)=u_2(x)y+xu_3(y)\; \forall x,y\in\AA\}.$$ 
By construction, $\ft(\AA)$ has three natural actions on $\AA$, which we denote $\AA_1, \AA_2,
\AA_3$. By a special case of the triality construction of Freudenthal's magic square 
\cite{LMdel}, 
we have a natural decomposition 
$$\fso(\AA)=Der\JA=\ft(\AA)\op\AA_1\op\AA_2\op\AA_3.$$
The Lie bracket can be explicitely described in this decomposition, but what we will
need again and again in the sequel is the explicit action on $\JA$, which is given by
the following formulas \cite{jacobson,harvey}. Let $u=(u_1,u_2,u_3)\in\ft(\AA)$ and 
$a_i\in\AA_i$. Then
 
\begin{eqnarray}
u\begin{pmatrix} r_1 & x_3 & x_2 \\ 
\overline{x_3} & r_2 & x_1 \\ \overline{x_2} & \overline{x_1} & r_3
\end{pmatrix} & = &
\begin{pmatrix} 0 & u_3(x_3) & u_2(x_2) \\ 
\overline{u_3(x_3)} & 0 & u_1(x_1) \\ \overline{u_2(x_2)} & \overline{u_1(x_1)} & 0
\end{pmatrix}, \\
 & & \nonumber \\
a_1\begin{pmatrix} r_1 & x_3 & x_2 \\ 
\overline{x_3} & r_2 & x_1 \\ \overline{x_2} & \overline{x_1} & r_3
\end{pmatrix} & = &
\begin{pmatrix} 0 & -\bar{a_1}\bar{x_2} & \bar{x_3}\bar{a_1} \\ 
-x_2a_1 & -2q(a_1,x_1) & (r_2-r_3)a_1 \\ a_1x_3 & (r_2-r_3)\bar{a_1} 
& 2q(a_1,x_1)\end{pmatrix}, \\
 & & \nonumber \\
a_2\begin{pmatrix} r_1 & x_3 & x_2 \\ 
\overline{x_3} & r_2 & x_1 \\ \overline{x_2} & \overline{x_1} & r_3
\end{pmatrix} & = &
\begin{pmatrix} 2q(a_2,x_2) & \bar{a_2}\bar{x_1} &(r_3-r_1)a_2   \\ 
x_1a_2 & 0 & -\bar{x_3}\bar{a_2} \\ (r_3-r_1)\bar{a_2} & -a_2x_3 
& -2q(a_2,x_2)\end{pmatrix}, \\
 & & \nonumber \\
a_3\begin{pmatrix} r_1 & x_3 & x_2 \\ 
\overline{x_3} & r_2 & x_1 \\ \overline{x_2} & \overline{x_1} & r_3
\end{pmatrix} & = &
\begin{pmatrix} 2q(a_3,x_3) & (r_2-r_1)a_3 & \bar{a_3}\bar{x_1} 
\\ (r_2-r_1)\bar{a_3} & -2q(a_3,x_3) & -\bar{x_2}\bar{a_3} \\ 
x_1a_3 & -a_3x_2 & 0\end{pmatrix}.
\end{eqnarray}
Recall that $q$ is the natural scalar product  on the complexified normed algebra $\AA$.
Using these formulas we can easily prove the Lemma. \medskip
 
\proof Since the action of $\SLA$ on the set of simply trisecant lines to $\DA$, we can
suppose that $L$ is the line of trace zero diagonal matrices, whose intersection with 
$\DA$ is the triple of points in $\PP\JA$ defined by the matrices $diag(0,1,-1), 
diag(1,0,-1)$ and $diag(1,-1,0)$. We then read off the previous formulas that the
subalgebra  of $\fso(\AA)$ consisting of operators that kill every diagonal (traceless)
matrices, is exactly $\ft(\AA)$. This implies the claim. \qed

\medskip
Note that $T^a_{w,v}=T^a_{w',v'}$ for any other couple
$(w',v')$ of rank 3 points of $\PP^2_{\e}$, such that the line 
$L' = \overline{w'v'}$ does not pass through a vertex $\e_i$ of $\Delta$. 
Indeed, the induced action 
of $T^a_{w,v} \subset \SLA$ on $\PJA$ must fix together with 
$w,v$ also all the lines (hence -- all the points) in the unique 
reduction plane  $\PP^2_{\e}$ passing through the line $L = \overline{wv}$.

\smallskip 
 {\sl Therefore, for any reduction plane $\PP^2_{\varepsilon} \subset \PJA$
the intersection 
$$
T^a_{\e} = \cap \{ \SOA_u : u \in \PP^2_{\varepsilon} - \Delta \}
$$
is, up to a finite group, 
 a copy of the triality group of the $\00$-column of the Freudenthal square. 
Moreover, 
$$
T^a_{\e} = T^a_{w,v} = \SOA_w \cap \SOA_v
$$
for any pair $w,v \in \PP^2_{\e}$, $w \not= v$, such that the line 
$L = \overline{wv}$ does not pass through a vertex $\e_i \in \AP2$ of the triangle 
$\Delta = \PP^2_{\e} \cap \DA$.}

\subsubsection{Lines, reduction planes and triality subspaces} 

Let $w,v \in \PP_{\e} - \Delta$ be as above, and let 
$\PP(V^a_w)$ and $\PP(V^a_v)$ be their isotropic subspaces 
of $\subset \PP(V^a)$. 
Then the intersection subspace 
$\PP(V^a_{\e}) = \PP(V^a_w) \cap \PP(V^a_v)=\cap \{ \PP(V^a_u) : u \in 
\PP^2_{\e} - \Delta \}$ in $\PP(V^a)$
is a projective representation space of their common triality 
subgroup $T^a_{w,v} = T^a_{\e}$.  

\begin{defi} {\bf (Triality subspaces)} 
For any reduction plane 
$\PP^2_{\e} = \langle \e_1;\e_2,\e_3\rangle \subset \PJA$ 
we call  the projective subspace $\PP(V^a_{\e}) \subset \PP(V^a)$ 
the triality subspace of $\PP^2_{\e}$.\end{defi}

\begin{prop} Let $a>1$, with the same notations as above.  
\begin{enumerate}
\item $Z^a_{w,v} := X^a_w \cap X^a_v = Z(w) \cap Z(v)$
is a copy of the triality variety $Z^a = \vee^3\PP^o, \times^3 \PP^1, D_4^{ad}$
from the $\underline{0}$-th column of the Freudenthal square. 
\item $Z^a_{w,v} = Z^a_{\e} := \cap \{ X^a_u: u \in \PP^2_{\e} - \Delta \}
                                            = \cap \{ Z(u) :  u \in \PP^2_{\e} - \Delta \}$.
\end{enumerate}
\end{prop} 

\proof
The second assertion follows from the previous discussion. Let us prove the 
first assertion. The case $a=2$ is easy, so we begin with $a=4$: here $w,v$ are 
generic symplectic forms on $\CC^6$, and $Z^4_{w,v}$ is the intersection in $G(3,6)$ 
of the lagrangian grassmannians $LG_w(3,6)$ and $LG_v(3,6)$. It is a classical fact that 
$w$ and $v$ can be simultaneously diagonalized, which means that we can find three
planes $P_1$, $P_2$, $P_3$ in general position in $\CC^6$, which are orthogonal 
with respect to both $w$ and $v$. It is then easy to see that a $3$-plane which
is isotropic with respect to $w$ and $v$ must be generated by three lines $l_i\subset 
P_i$. Hence the isomorphism $LG_w(3,6)\cap LG_v(3,6)\cong\times^3 \PP^1$. See \cite{ir}
for more details. 

Suppose now that $a=8$, and let again $L=\overline{wv}$ be the line of traceless
diagonal matrices in $\JO$. By the description of $\cE^8$ given in 
{\bf 2.5.2}, we have $Z(w)\cap Z(v)=E_6^{ad}\cap\PP\fso_8=D_4^{ad}$. \qed

\subsection{Triality varieties as zero-sets}

Let $w \in \PJA^o$, $a>1$; and denote by 
$\cF^a_w = \cE^a|_{X^a_w}$ the restriction of the homogeneous vector bundle
$\cE^a \rightarrow X^a$, to the isotropic subvariety  $X^a_w \subset X^a$. 
One can check case-by-case that $\cF^a_w$ is an irreducible homogeneous vector
bundle on $X^a_w$, and that its space of global sections is the polar hyperplane
$\JA_w$. Since $X^a_w$ is defined as the zero-locus of $w$, considered as a global 
section of $\cE^a$ on $X^a$, this isomorphism comes from the natural maps
$$\JA_w\simeq \JA/{\bf C} w\simeq H^0(X^a,\cE^a)/{\bf C}w\stackrel{res_{X^a_w}}{\lra}
 H^0(X^a_w,\cF_w^a).$$

For any $w \in \PJA^o$ 
the polar hyperplane $\PJA_w \subset \PJA$ to $w$ 
does not pass through $w$, and we can identify the base 
of the rational projection $p_w$ of $\PJA$ from $w$ 
with $\PJA_w$.
This identifies the lines $L \subset \PJA$ which pass through 
$w$ and the projective sections $w_L \in  \PP(X^a_w, \cF^a_w)$,
i.e. 
$$
p_w : \PJA \rightarrow \PJA_w = \PP(X^a_w, \cF^a_w), \quad
L \mapsto  w_L := p_w(L). $$
We shall denote by $Z(w_L) \subset X^a_w$ the zero-set of the 
projective section $w_L$ of $\cF^a_w$. 

For the fixed $w \in \PJA^o$,  let $U_w \subset G(2,\JA)$ 
be the subset of  all simply 3-secant lines to $\DA$ 
which pass through $w$. 
In  particular $p_w$ embeds $U_w$ in $\PJA_w$, 
and its complement 
$$
\Delta(\AA)_w = \PJA_w - U_w \subset \PJA_w
$$
is a hypersurface of degree $6$ in $\PJA_w$ 
 -- the {\it discriminant} hypersurface of $w$.

\smallskip

Let $L \in U_w$. 
By the preceding one can identify $L$ with a projective 
section 
$w_L \in U_w \subset \PJA_w = \PP H^0(X^a_w,\cF^a_w)$;
and we shall identify the zero-set  $Z(w_L) \subset X^a_w$ 
of $w_L$.  First, by Lemma \ref{3secant}, through the line 
$L \subset \PJA$ passes a unique simply 3-secant plane 
$\PP^2_{\e} = \PP^2_{\e(L)}$ to $X_a$; 
let $\PP^1_{\e(L)} \in Y_a^o$ be the proper 
$p_w$-image of $ \PP^2_{\e(L)}$ in $G(2,\JA_w)$. 
Second, by Proposition 2.9, for any $v \in \PJA^o$
the zero-set $Z(v) \subset X^a$ coincides with the 
isotropic subvariety $X^a_v \subset X^a$. 
Third, by Proposition 2.14, 
for any $v \in L^o =  L \cap \PJA^o$ the intersection 
$Z^a_L := X^a_w \cap X^a_v$
is a copy of the triality variety 
$Z^a$, embedded in the isotropy subvariety 
$X^a_w \subset X^a$; equivalently 
$$
Z^a_L = Z(w_L) =  \cap \{ Z(u) : u \in L^o \}
$$ 
is a copy of the triality variety 
$Z^a \subset X^a_w$. 

\subsection{The zero-set map}

We can identify the lines $L_t$ in $\PP^2_{\e(L)}$ passing through the
point $w$, with the elements of the line 
$\PP^1_{\e(L)} = p_w(\PP^2_{\e(L)}) \subset \PJA_w$.
Among these lines only three -- the lines 
$L_{\e_i} = \overline{w\e_i}, i = 1,2,3$ do not belong to $U_w$. 
Let $\overline{\e}_i = w_{L_{\e_i}}$ be their projections to 
$\PJA_w$. Then
$$
Z^a_L  = \cap \{ Z(w_{L_t}) : p_w(L_t) \in \PP^1_{\e(L)} - \overline{\e}(L) \}
$$ 
where $\overline{\e}(L)=\{\overline{\e}_1,\overline{\e}_2,\overline{\e}_3\}$. 

\medskip
Through any point 
$w_L \in \PP H^0(X^a_w,\cF^a_w)^o$
passes a {\it unique} 3-secant line $\PP^1_{\e} = \PP^1_{\e(L)}$ -- 
the proper $p_w$-image of the unique 3-secant plane 
$\PP^2_{\e(L)}$ which passes through the line $L$; 
and this plane is simply 3-secant to the Severi variety $X^a$.  
In other words, 
if $\Delta(\AA)_w \subset \PJA_w$ 
is the discriminant sextic, then the open subset 
$$
\PP H^0(X^a_w,\cF^a_w)^o  = \PJA_w^o = \PJA_w - \Delta(\AA)_w
$$
is swept out \underline{once} by the family of $3$-secant lines 
$\PP^1_{\e}$ to $\overline{X}^a_w$; and the unique line which
passes through a point of  $\PP H^0(X^a_w,\cF^a_w)^o$
belongs to the open subset $Y^o_a \subset Y_a$ of non-singular
reductions of $w$. 
One can regard such a 3-secant line $\PP^1_{\e}$ equivalently as a 
reduction $\varepsilon$ of $w$.  
By the preceding, all the points on such 3-secant 
line $\PP^1_{\e}$, except the 3 points of intersection of $\PP^1_{\e}$ 
with $\overline{X}_a$, are projective sections of $\cF^a_w$ with the 
same zero-set -- a triality variety $Z_{\e}$, defined uniquely 
by the reduction $\e$.  
Moreover, two 
different simply 3-secant lines $\PP^1_{\e}$ and $\PP^1_{\sigma}$
(equivalently -- two different non-singular reductions 
$\e = (\e_1,\e_2,\e_3)$ and $\s=(\sigma_1,\sigma_2,\sigma_3)$
of $w$) have different triality varieties $Z_{\e}$ and $Z_{\sigma}$. 
Indeed, two different reduction planes 
$\PP^2_{\e} = \langle\e_1,\e_2,\e_3\rangle$ and  
$\PP^2_{\sigma} = \langle\sigma_1,\sigma_2,\sigma_3\rangle$
in $\PJA$ define two different triality subgroups $T^a_{\e}$ and
$T^a_{\sigma}$  of $\SOA$, hence two different triality subspaces
$\PP(V^a_{\e})$ and $\PP(V^a_{\sigma})$, and two different triality 
subvarieties  $Z^a_{\e} \subset \PP(V^a_{\e})$ and 
$Z^a_{\sigma} \subset \PP(V^a_{\sigma})$. 

We collect all these observations in the following 

\begin{prop}
Let $a>1$, and let $w \in \PJA^o$. 
Let $p_w : \PJA \rightarrow \PJA_w$ be the rational 
projection from $w$ to the polar hyperplane section 
to $w$, and let $\overline{X}_a \subset \PJA_w$ be the 
isomorphic projection of the Severi variety $X_a \subset \PJA$ 
from $w$. 
Let $Y^o_a \subset Y^a$
be the open subset of non-singular reductions
$\e = (\e_1,\e_2,\e_3)$ of $w$; 
and denote by 
$$
\PP^2_{\e} = \langle\e_1,\e_2,\e_3\rangle \subset \PJA
\ \mbox{ and } \ \ 
\PP^1_{\e} =  p_w(\PP^2_{\e}) \subset \PJA_w
$$
correspondingly the plane  and the line of the 
non-singular reduction $\e \in Y_a^o$. 

Denote by ${\mathcal Z}(X^a_w)^o$ the (non closed) family of zero-sets 
$Z(v)$, for projective sections $v \in \PP H^0(X^a_w,\cF^a_w)^o$. 
Then any $Z(v) \in {\mathcal Z}(X^a_w)^o$ 
is a copy of the triality variety $Z^a$ in $X^a_w$; 
and there exists a 1:1-correspondence
$$
Y_a^o \cong {\mathcal Z}(X^a_w)^o, \  
\ \e \leftrightarrow Z_{\e}$$ 
described as follows:

{\bf 1.} 
If $\e \in Y_a^o$, then 
$Z_{\e} = \cap \{ Z(v): v \in {\PP}^1_{\e} - \overline{\e} \}$.

{\bf 2.}
If $Z \in {\mathcal Z}(X^a_w)^o$, then the closure of 
$$  \{ v \in \PP(H^0(X^a_w, \cF^a_w)): Z(v) = Z_{\varepsilon} \}
\ \subset \PP H^0(X^a_w, \cF^a_w) = \PJA_w  
$$ 
is a $3$-secant line $\PP^1_{\e}$ to $\overline{X}_w$, 
corresponding to a non-singular reduction $\e \in Y_a^o$ 
such that $Z = Z_{\e}$.
\end{prop} 

This suggests  to compactify $Y_a^0$ by embedding it into the Hilbert scheme
of $X^a_w$ and taking the closure, but we will not do that. 

We can also look at the sections $Z(v)$ of projective sections $v \in \PP H^0(X^a_w,\cF^a_w)$
that do not belong to the open subset $\PP H^0(X^a_w,\cF^a_w)$. An interesting ``degeneration'' 
of this kind occurs when $v$ is on a line joining, in some reduction plane $\PP^2_{\e}$, 
$w$ to one of the three vertices $\e_i$ of the simplex $\Delta =L_1+L_2+L_3$. 
We obtain the following picture:
$$\begin{array}{cccc}
a & 2 & 4 & 8 \\
Z_{\e} & \vee^3\PP^0 & \times^3\PP^1 & D_4^{ad} \\
Z(v) & \PP^1\vee \PP^0 & Q^3\times\PP^1 & B_4^{ad} 
\end{array}$$
Note that we have three vertices of $\Delta$, hence three varieties $Z(v_1), Z(v_2), Z(v_3)$
containing $Z_{\e}$. For $a=2$, we just have three points and we join two of them by a line.
For $a=4$ we have a product of three $\PP^1$'s, and we embed the product of two of them into
a three dimensional quadric. For $a=8$, $\ft(\OO)=\fso_8\subset\ft(\OO)\op\OO_i\cong\fso_9
\subset \ff_4$ for each $i=1,2,3$, so there are naturally three copies of $Spin_9$ in $F_4$ 
containing a given $Spin_8$, and there are three ways, inside $F_4^{ad}$, to embed the adjoint 
variety $D_4^{ad}=G_Q(2,8)$ in a copy of $B_4^{ad}=G_Q(2,9)$. Once again, triality leads the game. 

\vfill\pagebreak

\section{Geometry of the varieties of reductions}

In this section we make a detailed study of our varieties of reductions. 

\subsection{Varieties of reductions are linear sections of Grassmannians}

Our first result is that, as it is well-known for the Fano threefold $Y_1$, the 
varieties of reductions are linear sections of the ambient Grassmannians. But
there is a first surprise:

\begin{prop}\label{nontransverse} 
The variety of reductions $Y_a$ is a linear section of 
the Grassmannian $G(2,\JA_0)$, 
but non transverse for $a>1$, and not even of the expected dimension. 
\end{prop}

This will be proved in the next section. 

The linear section is defined as follows. Recall that the automorphism group 
${\rm Aut}\JA=\SOA_I$ preserves the quadratic form 
$Q(M)=\trace{(M^2)}$ on $\JA_0$.
At the infinitesimal level, this implies that the action of 
$\SOA_I$ on $\JA_0$  induces a map
from the Lie algebra $\fsoa$ of $\SOA$ to the space of skew-symmetric 
endomorphisms $\wedge^2\JA_0$. Explicitely, choose any orthonormal basis 
$X_i$ of $\JA_0$; then the map is 
$$u\in\fsoa\mapsto\sum_iX_i\wedge uX_i.$$
In particular, by Schur's lemma 
the wedge power $\wedge^2\JA_0$ contains a copy of $\fsoa$ as submodule.
It turns out that there is a very simple decomposition 
$$\wedge^2\JA_0=\fsoa\oplus U_a,$$
where the module $U_a$ is given as follows:
$$\begin{array}{ccccc}
a & 1 & 2 & 4 & 8 \\
\fsoa & \fsl_2 & \fsl_3 & \fsp_6 & \ff_4 \\
U_a & V_{6\om_1} & V_{3\om_1}\oplus V_{3\om_2}  & V_{\om_1+\om_3} & V_{\om_3}
\end{array}$$
Here we denoted by $V_{\om}$ the irreducible $\fsoa$-module of highest weight 
$\om$, and we used the same indexing of the weights as \cite{bou}. 

We can define the projection map $\pi :\wedge^2\JA_0\ra\fsoa$ by chosing two bases
$u_i$, $v_i$ of $\fsoa$ which are dual one to each other with respect to the Killing form. 
Then we can let, for $X,Y\in\JA_0$, 
$$\pi(X\wedge Y)=\sum_iQ(X,u_iY)v_i \in\fsoa.$$
It is quite clear that, the quadratic form $Q$ being $\fsoa$-invariant, this map
is well-defined, and equivariant. But by Schur's lemma there is only one such map, up to scalar, so 
$\pi$ must be the projection, up to scalar. This gives in particular a simple characterization 
of $U_a$, since it is precisely the kernel of $\pi$: it is the subspace of $\wedge^2\JA_0$ 
generated by the skew tensors $X\wedge Y$ such that $Q(X,uY)=0$ for all $u\in\fsoa$.

Note also that $U_a$ is always irreducible as a $\fsoa\times H_a$-module,
where $H_a$ is the finite group defined as follows (it is non trivial only for 
$a=2$, in which case $H_a=\ZZ_2$) : let $D_a$ be the Dynkin diagram 
of $\fsoa$, let $d_a\subset D_a$ be the set of nodes supporting the 
highest weight of $\JA_0$; then $H_a$ is the group of diagram 
automorphisms of $D_a$ preserving $d_a$. 

\smallskip The fact that $Y_a\subset \PP U_a$ can be seen as follows: if $\PP^2_{\e}$ 
is a simple trisecant plane to $X_a$ passing though $I$, and $\e_1, \e_2, \e_3$
denote the three intersection points, its representative in $\PP\JA_0\subset
\PP \JA$ is easily computed to be
$$\om_{\e}=\trace (\e_1)\e_2\wedge\e_3+\trace (\e_2)\e_3\wedge\e_1+
\trace (\e_3)\e_1\wedge\e_2.$$
Suppose that $\PP^2_{\e}$ is the plane of diagonal matrices. Then we immediately 
read off the formulas {\bf 2.7.1} (1--4)  that $Q(\e_i,u\e_j)=0$ for all $i,j$ and all $u\in\fsoa$. 
But the set of 
simple reduction planes is $\SOA$-homogeneous, so this is true in general. 
Thus each $\e_i\wedge\e_j$ is contained in $U_a$, and a fortiori $\om_{\e}$ also is. 

\subsection{Orbit structure}

The stabilizer $\SOA={\rm Aut}\JA$ of the identity element of $\JA$ acts on the 
variety of reductions. We prove that under this action, $Y_a$ only has a finite number 
of orbits.  More precisely:

\begin{prop}\label{orbits} 
The variety of reductions $Y_a$ is irreducible of dimension $3a$. It is the 
union of four $\SOA$-orbits of respective codimensions $0$, $1$, $2$ and $4$. 
(The codimension $4$ orbit is empty for $a=1$.)
\end{prop}

\proof We prove both Propositions simultaneously: we let 
$$\tilde{Y_a}:=G(2,\JA_0)\cap \PP U_a$$ 
and prove that is has four $\SOA$-orbits of codimensions $0$, $1$, $2$ and $4$. 
This will imply that $\tilde{Y_a}$ is irreducible. In particular it is equal to the closure
of its open orbit $Y_a^0$, the subset of reduction lines with three simple contacts with 
$\bar{X_a}$. Since $Y_a$ is precisely defined as the closure of this set of lines, it is 
equal to $\tilde{Y_a}$ and  Propositions \ref{nontransverse} and \ref{orbits} follow. 

\begin{lemm}\label{diag} Let $X\in\PP \JA_0-\Delta(\AA)$. Then there exists 
$g\in SO(\AA)_w$ such that $gX$ is a diagonal. \end{lemm}

\proof We first recall that the Cayley-Hamilton theorem holds
in $\JA$ \cite{chaput2}: any matrix $Y\in\JA$ satisfies  the identity
$$Y^3-\trace(Y)Y^2+Q'(Y)Y-\det(Y)I=0,$$ 
with $Q'(Y)=\frac{1}{2}(\trace(Y)^2-\trace(Y^2)).$
In particular, if $Y$ belongs to $\JA_0$, then $Y^3-\frac{1}{2}Q(Y)Y-\det(Y)I=0$. 
Moreover, the discriminant hypersurface can be defined, as usual, by the condition 
that the characteristic polynomial has a multiple root. 

Let $\a_1,\a_2,\a_3$ denote the roots of the characteristic polynomial of $X$. 
Since $X\notin \Delta(\AA)$, they are distinct. Let 
$$\begin{array}{ccc}
\pi_1=\frac{(X-\a_2 I)(X-\a_3 I)}{(\a_1-\a_2)(\a_1-\a_3)}, &
\pi_2=\frac{(X-\a_1 I)(X-\a_3 I)}{(\a_2-\a_1)(\a_2-\a_3)}, &
\pi_3=\frac{(X-\a_1 I)(X-\a_2 I)}{(\a_3-\a_1)(\a_3-\a_2)}.
\end{array}$$
A little computation shows that $\tr(\pi_i)=1$ and $\pi_i^2=\pi_i$. Moreover, 
$\pi_1+\pi_2+\pi_3=I$. In particular, $\langle\pi_1,\pi_2,\pi_3\rangle$ is a 
plane through $I$ with three simple contacts on $X_a$. 

But we know that $\SOA$ acts transitively on this set of planes, so that there
exists $g\in \SOA$ such that $g\langle\pi_1,\pi_2,\pi_3\rangle$ is the plane
of diagonal matrices. Since $X=\a_1\pi_1+\a_2\pi_2+\a_3\pi_3$, the matrix 
$gX$ is diagonal. \qed

\medskip
Let $l\in \tilde{Y_a}$ be a line which is not contained in the discriminant hypersurface
$\Delta (\AA)$. Choose a point $X\in l-\Delta (\AA)$. By the lemma, we can suppose that 
$X$ is diagonal. Then its diagonal coefficients are different, and the formulas 
{\bf 2.7.1} (1--4) imply that 
$$\fsoa X=\Bigg\{\begin{pmatrix} 0 & a_3 & a_2 \\ \bar{a_3} & 0 & a_1 \\
\bar{a_2} & \bar{a_1} & 0 \end{pmatrix}, \quad a_1,a_2,a_3\in\AA \Bigg\}.$$
Thus $\PP (\fsoa X)^{\perp}$ is the line of trace zero diagonal matrices (orthogonality
is taken with respect to the invariant quadratic form). 
This line is the projection of the plane of diagonal matrices, which is simply 
trisecant to the Severi variety $\overline{X_a}$. 

Therefore, the open subset of $\tilde{Y_a}$, of lines not contained in the
discriminant hypersurface, is an $SO(\AA)_w$-orbit isomorphic to the open set
of simply trisecant planes in $Y_a$. 

\smallskip
Let now $l\in \tilde{Y_a}$ be a line contained in the discriminant hypersurface
$\Delta(\AA)$. The plane generated by $l$ and the identity $I$ has a double contact 
at least with $X_a$, at some point 
$Z\in X_a^0$. Since $\SOA$ acts transitively on $X_a^0$ we may suppose that 
$$Z= \begin{pmatrix} 1 & i & 0 \\ i & -1 & 0 \\ 0 & 0 & 0 \end{pmatrix}, 
\quad T_ZX_a=\Bigg\{\begin{pmatrix} r-iu_0 & u & v \\ \bar{u} & r+iu_0 & iv \\ \bar{v} 
& i\bar{v} & 0 \end{pmatrix}, \quad r\in {\bf C}, u,v\in\AA \Bigg\}.$$
Here $u_0$ denotes the real part of $u$. 

Let us choose a tangent line  generated by a pair $(u,v)$, and consider the plane $P$ 
generated by this line and the identity. Observe that if this plane is in $\tilde{Y_a}$,
then $u$ must be real, i.e. equal to its real part $u_0$. Indeed, the projection of $P$ 
to $\PP \JA_0$ is the line joining $Z$ to the matrix 
$$Y=\begin{pmatrix} \frac{r}{3}-iu_0 & u & v \\ \bar{u} & \frac{r}{3}+iu_0 & iv 
\\ \bar{v} & i\bar{v} & -\frac{2r}{3} \end{pmatrix}.$$
For $s=(s_1,s_2,s_3)\in\ft(\AA)\subset\fsoa$, we have $Q(sZ,Y)=iq(s_3(1),u)$. 
Since $s_3(1)$, the image of the unit element $1\in\AA$, by the skew-symmetric endomorphism
$s_3\in\fso (\AA)$,  can be any imaginary vector in $\AA$, this scalar product is identically
zero if and only if $u$ is real. Changing $Y$ into $Y+u_0Z$ we can then suppose that 
$u=0$. We call the set of such tangent directions through $Z$
the {\em restricted set of tangents}.  

Now, a simple computation shows that $det(aI+bY+cZ)=c(c+rb)^2$. For $r\ne 0$, 
the intersection $P\cap \DA$ is the union of a tangent line to $X_a$ and a double 
non tangent line, both  through $Y$. The non tangent line cuts $X_a$ again, outside $X_a^0$, 
at the unique point 
$$X=\begin{pmatrix} q(v) & iq(v) & -rv \\ iq(v) & -q(v) & -irv 
\\ -r\bar{v} & -ir\bar{v} & r^2 \end{pmatrix}.$$
For $r=0$, $P\cap Sec(X_a)$ is a triple line 
through $Y$ tangent to $X_a$. If $q(v)\ne 0$, this line meets $X_a$ only at $Y$, 
but if $q(v)=0$ it is contained in $X_a$. 

This gives three cases, and we must check that we obtain correspondingly three 
$SO(\AA)$-orbits in $\tilde{Y_a}$, and no more. 

\begin{lemm} The isotropy group of $Z$ in $SO(\AA)$ acts on the restricted 
set of tangent directions through $Z$ with exactly three orbits, respectively 
of codimension  $0$, $1$ and $2$.\end{lemm}

\proof The isotropy subalgebra of (the line directed by) $Z$ in $\fsoa$ is 
$$\begin{array}{l}
Iso_Z(\fsoa)=\Big\{(s,a_1,a_2,a_3)\in\fsoa=\ft(\AA)\op\AA_1\op\AA_2\op\AA_3, \\
 \hspace{5cm}
a_1=ia_2, is_3(1)=2Im(a_3)\Big\}.
\end{array}$$
Its action on the restricted set of tangent directions is given by the formulas
$$\begin{array}{rcl}
(s,0,0,b)\begin{pmatrix} r & 0 & v \\ 0 & r & iv 
\\ \bar{v} & i\bar{v} & 0 \end{pmatrix} & = & 
iRe(b)\begin{pmatrix} 0 & 0 & v \\ 0 & 0 & iv 
\\ \bar{v} & i\bar{v} & 0 \end{pmatrix} +
\begin{pmatrix} 0 & 0 & s_2(v)-\frac{1}{2}s_3(1)v \\ 0 & 0 & is_1(v)+\frac{i}{2}s_3(1)v 
\\ * & * & 0 \end{pmatrix}, \\
 & & \\
(0,ia,a,0)\begin{pmatrix} r & 0 & v \\ 0 & r & iv 
\\ \bar{v} & i\bar{v} & 0 \end{pmatrix} & = & 
\begin{pmatrix} 0 & 0 & -ra \\ 0 & 0 & -ira 
\\ -r\bar{a} & -ir\bar{a} & 0 \end{pmatrix},
\end{array}$$
(with $is_3(1)=2Im(b)$, and the $*$ being the conjugates of the entries in symmetric position). 
On the second formula, we can already see that the isotropy 
group of $Z$ acts transitively on the set of restricted tangent directions for which 
$r\ne 0$. 

When $r=0$, we have to use the first formula, and study the rank of the map $\varphi_v : \ft(\AA)
\ra\AA$ sending $s$ to $s_2(v)-\frac{1}{2}s_3(1)v$. We claim that this map is surjective 
when $q(v)\ne 0$: this will ensure that the isotropy group of $Z$ acts transitively on the set 
of tangent directions for which $r=0$ and $q(v)\ne 0$. 

Recall from \cite{LMpop} that the triality algebra $\ft(\AA)$ is isomorphic to the direct 
sum of the derivation algebra $Der(\AA)$, with two copies of ${\rm Im}\AA$. Explicitely, the 
map sending $(D,u,v,w)\in Der(\AA)\op ({\rm Im\AA})^3$, with $u+v+w=0$, to the triple 
$s=(D+L_v-R_w,D+L_w-R_u,D+L_u-R_v)$, is an isomorphism onto $\ft(\AA)\subset\fso(\AA)^3$. 
(We denoted by $L_z$ and $R_z$ the operators of left and right multiplication by $z$ in 
$\AA$.)  We have 
$$s_2(v)-\frac{1}{2}s_3(1)v=Dv-\frac{1}{2}(tv+2vt),$$
so that the corank of $\varphi_v$ is equal to the corank of the endomorphism $\psi_v$ of 
$\AA$ defined by $\psi_v(t)=tv+2vt$. 

Suppose that $t\in Ker(\psi_v)$. Since $\AA$ is always alternative, the subalgebra 
generated by $t$ and $v$ is associative and we deduce that $2vtv=-tv^2=-4v^2t$. 
But since $v$ is imaginary, $v^2=-q(v)$, thus when $q(v)\ne 0$ we get $t=0$, as 
claimed. It follows that $\psi_v$ and $\varphi_v$ are surjective.

Now suppose that $q(v)=0$. Then the corresponding tangent direction is in fact the 
direction of a line which is contained in $X_a^0$. The
family of such lines is empty for $a=1$, and for $a=2$ it is the union of a projective
plane and its dual. For $a>2$, we know from \cite{LMhom}, Theorem 4.3, 
 that the family of lines in $X_a^0$  through the point $Z$
is irreducible, but splits into two orbits of the isotropy group, giving two types of lines 
which we called {\it general} and {\it special}, respectively. Already for dimensional
reasons we can see that the restricted tangent directions generate special lines only, 
hence that the isotropy group acts transitively on the set of tangent directions for
which $r=q(v)=0$. \qed

\medskip
We have thus obtained four orbits in $\tilde{Y_a}$, of codimension $0$, $1$, $2$ 
and $4$. It is clear from the proof that each orbit is in the closure of any 
other orbit of larger dimension. In particular, $\tilde{Y_a}$ is irreducible. 
This concludes the proof of Propositions \ref{nontransverse} and \ref{orbits}. \qed

\medskip Note that the identity $Y_a\cong \tilde{Y_a}$ 
implies the following characterization of lines belonging to $Y_a$, which we will use over and 
over in the sequel. 

\begin{coro}\label{quad} 
A line $\overline{XY}\subset\PP\JA_0$ defines a point of the variety of reductions
$Y_a\subset G(2,\JA_0)$,  if and only if 
$$Q(X,uY)=0\quad \forall u\in\fsoa.$$
\end{coro}

\medskip For future use we retain the following description of the $SO(\AA)_w$-orbits
in $\PP \JA_0$. We denote by ${\rm Tan}^0(X_a^0)\subset\PP\JA_0$ the union of the tangent 
lines to $X_a^0$ corresponding to the codimension one orbits of the isotropy groups of
the points of $X_a^0$, see the previous Lemma.  

\begin{prop}\label{J0orbits}
 The orbits of $SO(\AA)_w$ in $\PP \JA_0$ which are not contained in the 
discriminant hypersurface are the hypersurfaces $D_{[s,t]}(\AA)=\{X\in\PP \JA_0,\; 
6t\det(X)^2=sQ(X)^3\}$, where $[s,t]\in\PP^1-\{[1,9]\}$. 

The orbits of $SO(\AA)_w$ in $\PP \JA_0$ which are contained in the 
discriminant hypersurface $\Delta(\AA)=D_{[1,9]}(\AA)$ are: 
$$X_a^0, \quad \overline{X_a}-X_a^0,\quad  {\rm Tan}^0(X_a^0)-X_a^0,\quad \Delta(\AA)-\overline{X_a}
\cup {\rm Tan}(X_a^0).$$
\end{prop}

\proof The first assertion follows from Lemma \ref{diag}. To prove the second assertion, we first
recall that the action of $\SOA$ on $X_a$ (or $\overline{X_a}$) has exactly two orbits: the 
hyperplane section $X_a^0$ and its complement, see Lemma \ref{2orbits}. 
\smallskip
Let now $X\in\Delta(\AA)$
with $Q(X)\ne 0$. Then $Q(X)=6t^2$ and $\det (X)=-2t^3$ for a unique scalar $t$, 
and $(X-tI)^2(X+2tI)=0$. 
Let $Z=(X-tI)(X+2tI)$. If $Z=0$, then $X-tI$ is in $X_a$, thus $X$ belongs to $\overline{X_a}$. 
If $Z\ne 0$, $\trace (Z)=\trace(X^2)-6t^2=0$ and $Z^2=0$, so that $Z$ defines a point of 
$X_a^0$. Moreover, if $U=X+2tI$, we have $(U-3tI)Z=0$, hence $2UZ=6tZ=\trace (U)Z$, which means 
that $U$ belongs to  $T_ZX_a$. Since $X_a^0$ is $\SOA$-homogeneous, we may suppose that 
$$Z= \begin{pmatrix} 1 & i & 0 \\ i & -1 & 0 \\ 0 & 0 & 0 \end{pmatrix}, 
\qquad U=\begin{pmatrix} r-iu_0 & u & v \\ \bar{u} & r+iu_0 & iv \\ \bar{v} 
& i\bar{v} & 0 \end{pmatrix},$$
for some $r\in {\bf C}, u,v\in\AA$. Then we compute that the identity $Z=(X-tI)(X+2tI)=(U-3tI)U$
is equivalent to the relations $Im(u)=0$, $q(v)=0$, $r=3t$, $3tu_0=i$. Then we can write
$$X= \frac{i}{3t}Z+\begin{pmatrix} 0 & 0 & v \\ 0 & 0 & iv \\ \bar{v} 
& i\bar{v} & 0 \end{pmatrix}+t\begin{pmatrix} 1 & 0 & 0 \\ 0 & 1 & 0 \\ 0 & 0 & -2 \end{pmatrix},$$ 
where $q(v)=0$. If $v=0$, then again $X$ belongs to $\overline{X_a}$. If $v\ne 0$, by
the proof of Proposition \ref{orbits} the isotropy group of $Z$ acts transitively
on the special lines of $X_a^0$ passing through $Z$. Looking more carefully at the explicit
action, we can see that in fact, it acts transitively on the cone of $X_a^0$ generated by these
lines (minus the vertex $Z$, of course). We conclude that 
$\SOA$ acts transitively on $\Delta(\AA)-Tan^0(X_a^0)\cup\overline{X_a}$. 

\smallskip
Let us consider now some  $X\in\Delta(\AA)$ such that $Q(X)=\det (X)=0$. In particular, 
$X^3=0$, and $X$ belongs to $X_a^0$ if $X^2=0$. Suppose this is not the case. 
Then $Z=X^2$ defines 
a point of $X_a^0$, and $U=X$ belongs to $T_ZX_a^0$. We can chose $Z$ to be the same matrix as in 
the previous case, as well as $U$, but with $r=0$ since $\trace (U)=\trace (X)=0$. 
The equation $Z=U^2$ gives the relations $q(v)=1$ and $Im(u)v=0$, hence $Im(u)=0$ since
$v$ is invertible. Again, we check that the isotropy 
group of $Z$ acts transitively on the (pointed) cone generated by the codimension one orbit
of restricted tangent directions through $Z$. We conclude that 
$\SOA$ acts transitively on ${\rm Tan}^0(X_a^0)-X_a^0$, and the proof is complete. \qed

\medskip
Explicit representatives of the four orbits in $\Delta(\AA)$ are, respectively:
$$\begin{pmatrix} 1 & i & 0 \\ i & -1 & 0 \\ 0 & 0 & 0 \end{pmatrix}, \quad 
\begin{pmatrix} 1 & 0 & 0 \\ 0 & 1 & 0 \\ 0 & 0 & -2 \end{pmatrix}, \quad 
\begin{pmatrix} 1 & i & 1 \\ i & -1 & i \\ 1 & i & 0 \end{pmatrix}, \quad 
\begin{pmatrix} 1 & 1+iI & 0 \\ 1-iI & 1 & 0 \\ 0 & 0 & -2 \end{pmatrix}.$$
\medskip

\subsection{Geometric description of the orbits}
We denote by $Y_a^i$ the orbit of codimension $i$ in $Y_a$. 

\medskip\noindent {\sl a. The open orbit}. 

\begin{prop}\label{open}
The open $\SOA$-orbit in $Y_a$ is $Y_a^0\simeq \SOA/T_a$, with 
$Lie(T_a)=\ft(\AA)$. In particular, it is an affine variety. 
\end{prop}

\proof A point in the open orbit $Y_a^0$ is given by the line of 
traceless diagonal matrices in $\PP \JA_0$. One can read 
off the explicit action of $\fsoa$ on $\JA_0$ that the 
stabilizer of this line is $\ft(\AA)$, which implies the first
assertion. Since $\ft(\AA)$ is reductive, the second assertion  
follows from a theorem of Matsushima, following which the 
quotient of a reductive group by a reductive subgroup is affine \cite{matsu}. 
\qed 

\medskip
The open orbit $Y_a^0$ consists in planes in $\PP \JA$ having
three simple contacts with $X_a$. Let $\langle \e_1, \e_2, \e_3\rangle$
be such a plane, and suppose that $\e_1+\e_2+\e_3=I$. By Corollary 2.11, 
this triple is projectively equivalent to the triple of rank one diagonal 
matrices, and since $\SOA={\rm Aut}\JA$ this implies that $\e_1\e_2=\e_1\e_3=\e_2\e_3=0$.
Multiplying the previous identity by $\e_i$, we thus get $\e_i^2=\e_i$, 
hence $\trace(\e_i)=1$  by Lemma 2.2. The triple $ \e_1, \e_2, \e_3$ is therefore 
what algebraists call a {\it Pierce decomposition} of the Jordan algebra $\JA$ \cite{jacobson}. 

These observations lead to the slightly more precise statement:

\begin{prop} \label{s3}
The generic isotropy group $T_a$ is the semi-direct product of the triality 
group $T(\AA)$ with the symmetric group  ${\mathfrak S}_4$. 
\end{prop}

\proof The generic isotropy group $T_a$ is the stabilizer of the line of traceless
diagonal matrices in $\JA_0$, or equivalently to the plane of diagonal matrices
in $\JA$, which is generated by the three diagonal idempotents 
$\e_1,\e_2,\e_3$. 
These three idempotents are permuted by the elements of $T_a$, giving a 
morphism
$\nu :T_a\ra {\mathfrak S}_3$. Note that ${\mathfrak S}_3$ 
is the quotient of ${\mathfrak S}_4$ 
by the normal subgroup generated by the permutations which are products 
of two disjoint transposition. This subgroup is a copy of $\ZZ_2\times\ZZ_2$.  
We must therefore prove that the morphism $\nu$ is 
surjective, and that its kernel $T_a^0$ coincides with the 
semi-direct product of $T(\AA)$ with $\ZZ_2\times\ZZ_2$. 

To prove the surjectivity of $\nu$, we define two endomorphisms $\s_1$ and 
$\s_2$ of $\JA$ by the formulas
$$\begin{array}{rcl}
\s_1 \begin{pmatrix} 
c_1            & x_3            & x_2  \\ 
\bar{x}_3 & c_2            &  x_1            \\ 
\bar{x}_2           & \bar{x}_1 &  c_3            
\end{pmatrix}=
\begin{pmatrix} 
c_2          & \bar{x}_3           & x_2  \\ 
x_3 & c_1            &  x_1            \\ 
\bar{x}_2        & \bar{x}_1 &  c_3            
\end{pmatrix}, \\
 & & \\
\s_2 \begin{pmatrix} 
c_1            & x_3            & x_2  \\ 
\bar{x}_3 & c_2            &  x_1            \\ 
\bar{x}_2           & \bar{x}_1 &  c_3            
\end{pmatrix}=
\begin{pmatrix} 
c_1          &   x_3         & x_2  \\ 
\bar{x}_3 & c_3            &  \bar{x}_1            \\ 
\bar{x}_2        & x_1 &  c_2            
\end{pmatrix}.
\end{array}$$
It is easy to check that these endomorphisms are in fact automorphisms 
of the Jordan algebra $\JA$. Moreover, they belong to $T_a$, and their
images by $\nu$ are the two simple generators of ${\mathfrak S}_3$, proving that 
$\nu$ is surjective. 

Let now $t\in T_a^0$, so that $t$ fixes each $\e_i$. Being an automorphism of $\JA$,
$t$ also preserves the subspace $\e_i\JA \e_j$, for $1\le i,j\le 3$. 
But this is 
the space of matrices whose entries are zero except possibly that on the $i$-th line 
and $j$-th column, and symmetrically that on the $j$-th line 
and $i$-th column. Therefore, there exists scalars $\zeta_1, \zeta_2, \zeta_3$, 
and endomorphisms $\tau_1, \tau_2, \tau_3$ of $\AA$, such that 
$$t \begin{pmatrix} 
c_1            & x_3            & \bar{x}_2  \\ 
\bar{x}_3 & c_2            &  x_1            \\ 
x_2            & \bar{x}_1 &  c_3            
\end{pmatrix}=
\begin{pmatrix} 
\zeta_1c_1         & \tau_3(x_3)   &  \tau_2(x_2)  \\ 
\overline{\tau_3(x_3)} & \zeta_2c_2           & \tau_1(x_1)             \\ 
\overline{\tau_2(x_2)}     & \overline{\tau_1(x_1)} &  \zeta_3c_3            
\end{pmatrix}.$$
But $t(I)=I$, hence $\zeta_1=\zeta_2=\zeta_3=1$. Moreover, a straightforward
computation shows that this is an automorphism of $\JA$ if and only if 
$$\tau_2(xy)=\tau_3(x)\tau_1(y)\qquad \forall x,y\in\AA.$$
This is precisely the definition of $T(\AA)$, except that we don't ask the 
$\tau_i$ to belong to $SO(\AA)$. They will automatically belong to the 
orthogonal group $O(\AA)$, but the sign ambiguity explains the 
appearance of the $\ZZ_2\times\ZZ_2$ factor.  This
concludes the proof. (The heart of this argument can be found in 
\cite{jacobson}, see also \cite{sk}). \qed

\medskip\noindent {\sl b. The orbit of codimension one}. 

\smallskip
The codimension one orbit $Y_a^1$ is made of planes in $\PP \JA$ 
having a simple contact with $X_a$ outside $X_a^0$, and a double contact
on $X_a^0$. This implies the existence of two fibrations $p$ and $p_0$, 
the first one over $X_a-X_a^0$,
the second one over $X_a^0$.

\begin{lemm} {\bf (Point-line polarity in $X_a=\AA\PP^2$)}.  
Let $Z$ be a point of $X_a$.
The intersection $X_a\cap (T_ZX_a)^{\perp}$ is an $a$-dimensional
quadric $Q^a_Z$, an $\AA$-line in $X_a$. This quadric contains $Z$ if 
and only if $Z$ belongs to $X_a^0$. 

A point $Y$ belongs to $Q^a_Z$ if and only if there exists a reduction plane 
$P\in Y_a$ passing through $Y$ and $Z$. In particular, 
$$Y\in Q^a_Z \iff Z\in Q^a_Y.$$
\end{lemm}

\proof If $Z\notin X_a^0$, the line $\overline{ZI}$ meets the determinant
hypersurface $\DA$ at a unique point $M\notin X_a$, and the set of secant 
(or tangents) lines to $X_a$
passing through $M$ cuts a smooth $a$-dimensional quadric $\tilde{Q}^a_Z$ on 
$X_a$ (the entry-locus of $M$, see e.g. \cite{zak}). A reduction plane $P\in Y_a$
through $Z$ is then generated, with $Z$ an $I$, by a point $Y$ of that quadric. 
A simple computation shows that actually, $\tilde{Q}^a_Z=Q^a_Z$, and our claim 
follows for all points outside $X_a^0$. 

To conclude the proof, we check than also for $Z\in X_a^0$, the intersection 
$X_a\cap (T_ZX_a)^{\perp}$ is a smooth quadric: this is a straightforward
computation. Then the  last assertion of the Lemma follows by continuity.\qed 

\medskip\noindent {\it Remark}. If $Z\notin X_a^0$, one can define 
the quadric $Q^a_Z$ as the set of points $Y\in X_a$ such that $YZ=0$. 
Nevertheless, for $Z\in X_a^0$ this condition defines a larger set than
$Q_Z^a$.

\medskip
This Lemma allows a simple description of the fibers of the two projections 
$p$ and $p_0$. Indeed, one easily checks that if $Z\in X^0_a$ and $X\in X_a-X_a^0$, 
$$p_0^{-1}(Z)=Q^a_Z-Q^a_Z\cap X_a^0, \qquad {\rm and}\quad 
p^{-1}(X)=Q^a_X\cap X_a^0.$$
We'll see below that $Q^a_Z\cap X_a^0$ is a singular section of the quadric $Q^a_Z$, 
with a unique singularity at $Z$, so that $p_0^{-1}(Z)\simeq \CC^a$.

An explicit representative of $Y_a^1$ is the line generated by 
$$Z=\begin{pmatrix} 1 & i & 0 \\ i & -1 & 0 \\
0 & 0 & 0 \end{pmatrix}
\quad {\rm and}\quad 
Y=\begin{pmatrix} 1 & 0 & 0 \\ 0 & 1 & 0 \\
0 & 0 & -2 \end{pmatrix}.$$

\medskip\noindent {\sl c. The orbit of codimension two}. 

\smallskip
The elements of the codimension two orbit $Y_a^2$ are the planes in $\PP \JA$ 
with a triple contact with $X_a$ on $X_a^0$. This defines a fibration $p$ over
$X_a^0$, and by the proof of Proposition \ref{orbits},  the fiber 
of $p$ over $Z$ is the set of tangent directions generated by matrices of the form 
$$Y=\begin{pmatrix} 0 & 0 & v \\ 0 & 0 & iv \\ \bar{v} 
& i\bar{v} & 0 \end{pmatrix}\qquad {\rm if}\quad 
Z= \begin{pmatrix} 1 & i & 0 \\ i & -1 & 0 \\ 0 & 0 & 0 \end{pmatrix},$$
with $q(v)\ne 0$. An easy computation shows that 
$$Q^a_Z=\Bigg\{\begin{pmatrix} r_1 & ir_1 & x_2 \\ ir_1 & -r_1 & ix_2 \\ \bar{x_2} 
& i\bar{x_2} & r_3 \end{pmatrix}, \;\;r_1,r_3\in {\bf C}, x_2\in\AA, \; r_1r_3=q(x_2)\Bigg\}.$$
This is a smooth quadric which is tangent to $X_a^0$ at $Z$, hence $Q_Z^a\cap X_a^0$ 
is a quadratic cone with vertex $Z$, which is its unique singular point. The fiber 
$p^{-1}(Z)$ can then be described as the set of lines through $Z$ in the linear 
subspace of $\PP\JA_0$ spanned by $Q_Z^a\cap X_a^0$, which are not contained in that 
cone. This shows that $p^{-1}(Z)$ is the complement of a smooth quadric hypersurface 
in a $\PP^{a-1}$. 

An explicit representative of $Y_a^2$ is the line generated by 
$$Z=\begin{pmatrix} 1 & i & 0 \\ i & -1 & 0 \\
0 & 0 & 0 \end{pmatrix}
\quad {\rm and}\quad 
Y=\begin{pmatrix} 0 & 0 & 1 \\ 0 & 0 & i \\
1 & i & 0 \end{pmatrix}.$$

\medskip\noindent {\sl d. The closed orbit}. 
\smallskip

Finally, the codimension four orbit $Y_a^4$ is made of 
planes containing a line of $X_a^0$, and such planes are completely 
determined by their corresponding line, which we noticed to be special 
when $a>2$. With the convention that any line on $\times^2\PP^2$ is special,
we get:

\begin{prop} The closed orbit $Y_a^4$ in $Y_a$ is isomorphic to the orbit of 
special lines in $X_a^0$. \end{prop}

Specifically, we have
$$Y_1^4=\emptyset, \quad Y_2^4=\PP^2\sqcup{\check \PP}^2, \quad Y_4^4=F_{\om}(1,3;6)=Sp_6/P_{1,3},
\quad Y_8^4=F_4/P_3.$$

An explicit representative of $Y_a^4$ is the line generated by 
$$Z=\begin{pmatrix} 1 & i & 0 \\ i & -1 & 0 \\ 0 & 0 & 0 \end{pmatrix}
\quad {\rm and}\quad 
Y=\begin{pmatrix} 0 & 0 & 1-iI \\ 0 & 0 & i+I \\
1+iI & i-I & 0 \end{pmatrix}.$$

\subsection{Smoothness}
We have seen in Proposition \ref{nontransverse} that the varieties of reductions, for $a>1$, 
are non transverse linear sections of their ambient Grassmannians. The following Theorem
is therefore rather surprising. 

\begin{theo}\label{smoothness} The varieties of reductions $Y_a$ are smooth.\end{theo}

\proof This is already known for $a=1$. For $a>1$, we check that a point of the 
codimension $4$ orbit is smooth, which is enough to prove the theorem. We have 
just seen that 
a point of $Y_a^4$ is the line generated by 
$$Z=\begin{pmatrix} 1 & i & 0 \\ i & -1 & 0 \\ 0 & 0 & 0 \end{pmatrix}
\quad {\rm and}\quad 
Y=\begin{pmatrix} 0 & 0 & 1-iI \\ 0 & 0 & i+I \\
1+iI & i-I & 0 \end{pmatrix}.$$
We choose a basis $e_1=1+iI, e_2,\ldots ,e_a$ of $\AA$. Then we can complete these two matrices 
$Z,Y$ into a basis of $\JA_0$, 
$$\begin{array}{l}
Z=\begin{pmatrix} 1 & i & 0 \\ i & -1 & 0 \\ 0 & 0 & 0 \end{pmatrix}, \quad
Y_j^+=\begin{pmatrix} 0 & 0 & e_j \\ 0  & 0 & ie_j \\ \bar{e_j} & i\bar{e_j} & 0 \end{pmatrix}, 
\quad  
Y_j^-=\begin{pmatrix} 0 & 0 & e_j \\ 0  & 0 & -ie_j \\ \bar{e_j} & -i\bar{e_j} & 0 \end{pmatrix}, \\    \\
\hspace{1.9cm}
X_j=\begin{pmatrix} 0 & e_j & 0 \\ \bar{e_j} & 0 & 0 \\ 0 & 0 & 0 \end{pmatrix}, \qquad
W=\begin{pmatrix} 0 & 0 & 0 \\ 0  & 1 & 0 \\ 0 & 0 & -1 \end{pmatrix}. 
\end{array} $$ 
This provides us with a system of local coordinates on $G(2,\JA_0)$ around the 
line $\overline{ZY}$ : a line in a certain neighbourhood of $\overline{ZY}$ has a unique basis 
$A,B$ of the form
$$\begin{array}{rcl} 
A & = & Z
+\sum_{j>1}a_j^+Y_j^++\sum_{j\ge 1}a_j^-Y_j^-+\sum_{j\ge 1}r_jX_j+uW, \\
B & = & Y
+\sum_{j>1}b_j^+Y_j^++\sum_{j\ge 1}b_j^-Y_j^-+\sum_{j\ge 1}s_jX_j+vW.
\end{array}$$
By Lemma \ref{quad}, 
such a line belongs to $Y_a$ if and only if $Q(A,uB)=0$ for all $u\in\fsoa$. If we 
write $A=Z+\delta A$ and $B=Y+\delta B$, we get the infinitesimal equations 
$$Q(uZ,\delta B)=Q(uY,\delta A) \qquad \forall u\in\fsoa.$$
Using the explicit action of $m_i\in\AA_i\subset\fsoa$, we obtain the following 
three sets of $a$ equations,
$$\begin{array}{l}
2i\sum_{j\ge 1}q(m_1,e_j)b_j^- \quad = \quad 
  \sum_{j\ge 1}q(e_1\bar{m_1},e_j)r_j+4iq(m_1,e_1)u, \\
2i\sum_{j\ge 1}q(m_2,e_j)b_j^-  \quad = \quad  -\sum_{j\ge 1}iq(m_2\bar{e_1},e_j)r_j
+\sum_{j>1}iq(m_2,e_j)a_j^+ \\
  \hspace{5cm} +\sum_{j\ge 1}q(m_2,e_j)a_j^-+2iq(m_2,e_1)u, \\
-2\sum_{j\ge 1}q(m_3,e_j)s_j-2iq(m_3,1)v \quad = \quad \sum_{j>1}iq(m_3\bar{e_1},e_j)a_j^+ \\ \hspace{8cm}
+\sum_{j\ge 1}q(m_3\bar{e_1},e_j)a_j^-.
\end{array}$$
The first set of equations gives the $b_j^-$ in terms of $u$ and the $r_j$, because
the coefficient of $b_j^-$ is  $q(m_1,e_j)$, and $m_1$ can be chosen arbitrarily . Then the
second set of equations gives the $a_j^-$ in terms of $u$, the $r_j$ and the $a_j^+$. 
Finally, the third set of equations gives the $s_j$ in terms of $u,v$,  the $r_j$ and 
the $a_j^+$. 

This proves that the Zariski tangent space of $Y_a$ at $\overline{ZY}$ has codimension at least $3a$ 
in that of $G(2,\JA_0)$. But $Y_a$ has dimension $3a$ and $G(2,\JA_0)$ has dimension
$6a$, so $\overline{ZY}$ must be a smooth point of $Y_a$. \qed

\subsection{Linear spaces in the varieties of reductions}
In a Grassmannian $G(2,n)$ of projective lines, there are two types 
of linear spaces. Those {\it of the first kind} are made of the lines containing
a fixed point and contained in a fixed subspace. Those {\it of the second 
kind} have dimension two only; they are made of the lines contained in a fixed plane. 

By Lemma \ref{quad}, the maximal linear spaces of the first kind
that are contained in $Y_a$  are defined as follows: take some point
$X\in \PP \JA_0$ and consider the space of lines $L$ such that 
$$X\in  L\subset \PP(\fsoa X)^{\perp}.$$

\begin{prop}\label{subspaces} The space $(\fsoa X)^{\perp}$ has dimension $a+2$ 
if $X$ belongs to the projection $\overline{X_a}$ of the Severi 
variety $X_a$, and dimension $2$ otherwise.
\end{prop}

\proof The dimension of $\fsoa X$ is the dimension of the $SO(\AA)_w$-orbit of 
$X$ in $\JA_0$ (not in the projectivisation $\PP \JA_0$ !). We consider the different cases
obtained in Proposition \ref{J0orbits}. 

If $X\notin\Delta (\AA)$, its orbit is, by Lemma \ref{diag}, 
 the set of matrices with the same characteristic polynomial:
its codimension is two, hence $(\fsoa X)^{\perp}$ has dimension two. 

If $X$ belongs to (the cone over) the open orbit in $\Delta(\AA)$, again its orbit in $\JA_0$
depends on its characteristic polynomial, so it must be of  codimension two and again
$(\fsoa X)^{\perp}$ has dimension two. 

The (pointed) cone over ${\rm Tan}^0(X_a^0)-X_a^0$ either is an $\SOA$-orbit, or the
union of a one dimensional family of codimension one orbits. An explicit computation 
of the infinitesimal action shows that we are in fact in the first situation. Indeed, 
By Proposition 3.6 and its comment we can let
$$ X=\begin{pmatrix} 1 & i & 1 \\ i & -1 & i \\ 1 & i & 0 \end{pmatrix},$$
and determine its centralizer using formulas (1--4) in {\bf 2.7.1}. If $(u,a_1,a_2,a_3)\in
\fsoa$ annihilates $X$, looking at the diagonal coefficients we first see that $a_1$, $a_2$ and
$a_3$ must be imaginary. The non diagonal coefficients then give the equations
$$\begin{array}{rcl}
u_1(1) & = & -ia_1-a_2-ia_3, \\
u_2(1) & = & -a_1+a_2+ia_3, \\
u_1(1) & = & ia_1+a_2-2ia_3.
\end{array}$$
The matrix formed by the coefficients of $a_1,a_2,a_3$ is easiliy seen to be invertible.
We conclude that $a_1,a_2$ and $a_3$ are uniquely determined by $u$, which can
be arbitrary. Thus the stabilizer of $X$ has codimension $3a$ in $SO(\AA)$, which implies
that the orbit of $X$ has dimension $3a$, which is also the dimension of the cone over 
${\rm Tan}^0(X_a^0)$. This proves our claim. 

Hence if $X$ belongs ${\rm Tan}^0(X_a^0)-X_a^0$, since it is a codimension two orbit 
in $\JA_0$,  $(\fsoa X)^{\perp}$ has dimension two again. 

If $X$ belongs to the cone over $\overline{X_a}-X_a^0$, its $\SOA$-orbit 
depends on its characteristic polynomial, so its orbit has dimension $2a$. 
Finally, the (pointed) cone over $X_a^0$ is a full $\SOA$-orbit of dimension $2a$. 
Thus if $X$ belongs to the cone over $\overline{X_a}$, the dimension of $(\fsoa X)^{\perp}$ 
is $a+2$, independently of the fact that $X$ belongs to $X_a^0$ or not. \qed 

\begin{coro} The variety $Y_a\subset G(2,\JA_0)$ does not contain
any plane of the second kind. \end{coro}

\proof A plane of the second kind in $G(2,\JA_0)$ is a space of lines contained
in the projectivization of some three-dimensional subspace $K$ of $\JA_0$. 
Let such a plane be contained in $Y_a$, and consider a point of that plane, which 
represents a line $\overline{xy}$ in $\PP\JA_0$.
Then we must have $K\subset (\fsoa x)^{\perp}\cap (\fsoa y)^{\perp}$. In particular 
the line $\overline{xy}\subset \PP \JA_0$ must be contained in $\overline{X_a}$:
otherwise, if $p=\a x+\b y\notin \overline{X_a}$, then $(\fsoa x)^{\perp}\cap (\fsoa y)^{\perp}
\subset (\fsoa p)^{\perp}$, which is two-dimensional. 

Suppose that $x\in \overline{X_a}-X_a^0$. Since this space is homogeneous, we can let 
$$x=\begin{pmatrix} 1 & 0 & 0 \\ 0 & 1 & 0 \\ 0 & 0 & -2 \end{pmatrix}, \quad 
{\rm thus}\quad (\fsoa x)^{\perp}=\Bigg\{
\begin{pmatrix} r_1 & x_3 & 0 \\ \bar{x_3} & r_2 & 0 \\ 0 & 0 & r_3 \end{pmatrix}, 
r_1+r_2+r_3=0\Bigg\}.$$
A straightforward computation shows that 
$\PP (\fsoa x)^{\perp}\cap\overline{X_a}$ is
$$\Bigg\{
\begin{pmatrix} r_1 & x_3 & 0 \\ \bar{x_3} & r_2 & 0 \\ 0 & 0 & r_3 \end{pmatrix}, \; 
r_1+r_2+r_3=0, q(x_3)=(r_1-r_3)(r_2-r_3)\Bigg\}\cup\{x\}.$$
In particular, $x$ is an isolated point of that intersection, which can therefore contain 
no line through $x$. 

Suppose now that $x\in X_a^0$, and we can even suppose that the whole line $\overline{xy}$
is contained in $X_a^0$. By the proof of Proposition \ref{orbits}, we can suppose that 
$$x=\begin{pmatrix} 1 & i & 0 \\ i & -1 & 0 \\ 0 & 0 & 0 \end{pmatrix}, \qquad 
y=\begin{pmatrix} 0 & 0 & v \\ 0 & 0 & iv \\ \bar{v} & i\bar{v} & 0 \end{pmatrix},$$
where $q(v)=0$. Another computation, which we leave to the reader,  
shows that again $(\fsoa x)^{\perp}\cap (\fsoa y)^{\perp}$ is only two-dimensional. \qed 

\begin{coro} The maximal linear spaces in $Y_a$ are $\PP^a$'s
parametrized by  $X_a$. \end{coro}

\medskip Recall (Lemma \ref{2orbits}) that $\SOA$ has exactly two orbits inside $X_a\simeq 
\overline{X_a}$:
the closed orbit $X_a^0$, which is the hyperplane section of $X_a$ by
$\PP \JA_0$, and its complement. Correspondingly, there are two types
of $\PP^a$'s inside $Y_a$ : {\em special} ones, for $x\in X_a^0$, and 
{\em general} ones, for $x\notin X_a^0$.

\begin{prop} The numbers of general and special $\PP^a$'s through a point 
of the codimension $i$ orbit $Y_a^i$ of $Y_a$ is given as follows:
$$\begin{array}{lcc}
\hspace{1cm} & general & special \\
Y_a^0 & 3 & 0 \\
Y_a^1 & 1 & 1 \\
Y_a^2 & 0 & 1 \\
Y_a^4 & 0 & \infty^1 \end{array}$$
\end{prop}

\medskip Note that the fact that there are exactly three $\PP^a$'s 
through a point of the open orbit is a genuine geometric 
manifestation of triality !
Indeed, we know that the tangent space to a point of the open 
orbit is equal to $\fsoa/\ft(\AA)=\AA_1\op\AA_2\op\AA_3$ as a module
over the stabilizer Lie algebra $\ft(\AA)$. The three copies of $\AA$ 
correspond to the directions of the three $\PP^a$'s. Moreover, by Proposition \ref{s3}
the isotropy group of a generic point contains a copy of ${\mathfrak S}_3$, which
permutes these three spaces. \medskip

\proof A point of $Y_a^0$ is given by the line of diagonal matrices in 
$\PP \JA_0$. If this point is contained in a maximal linear subspace $\PP^a_x$
of $Y_a$, then $x$ is diagonal and belongs to $\overline{X_a}$. There are exactly
three such matrices (up to scalar), 
$$\begin{pmatrix} 1 & 0 & 0 \\ 0 & 1 & 0 \\ 0 & 0 & -2 \end{pmatrix}, \qquad 
\begin{pmatrix} 1 & 0 & 0 \\ 0 & -2 & 0 \\ 0 & 0 & 1 \end{pmatrix}, \qquad 
{\rm and} \quad 
\begin{pmatrix} -2 & 0 & 0 \\ 0 & 1 & 0 \\ 0 & 0 & 1 \end{pmatrix},$$
the projections of the three diagonal matrices of rank one. Thus a point of 
$Y_a^0$ belongs to exactly three lines, and they are all general. 

A point of $Y_a^1$ is given by the line in $\JA_0$ generated by 
$$Z=\begin{pmatrix} 1 & i & 0 \\ i & -1 & 0 \\
0 & 0 & 0 \end{pmatrix}
\quad {\rm and}\quad 
Y=\begin{pmatrix} 1 & 0 & 0 \\ 0 & 1 & 0 \\
0 & 0 & -2 \end{pmatrix}.$$
This line meets $\overline{X_a}$ only at $Z$ and $Y$. Since $Y$ belongs to $X_a^0$ and $Z$ 
does not, this implies that a point of $Y_a^1$ belongs exactly to one special $\PP^a$  and 
one general $\PP^a$ of $Y_a$. 

A point of $Y_a^2$ is given by the line generated by 
$$Z=\begin{pmatrix} 1 & i & 0 \\ i & -1 & 0 \\
0 & 0 & 0 \end{pmatrix}
\quad {\rm and}\quad 
Y=\begin{pmatrix} 0 & 0 & 1 \\ 0 & 0 & i \\
1 & i & 0 \end{pmatrix}.$$
This is a tangent line to $X_a^0$ at $Z$, and it meets $\overline{X_a}$ only at 
this point.  Thus a point in $Y_a^2$ belongs to a unique maximal linear subspace
of $Y_a$, which is special. 

Finally, a point  of $Y_a^4$ is given by a special line in $X_a^0$, and each point of this 
line defines a special maximal linear subspace of $Y_a$ through that point 
of $Y_a^4$. \qed

\subsection{Varieties of reductions are rational Fano manifolds}

Let $Z_a$ denote the space of incident points and lines
$p\subset l$, with $p\in\PP \JA_0$ and $l\in Y_a$. 

\begin{prop}\label{blowup}
There is a commutative diagram 
$$\xymatrix@1{  & Z_a=\PP_{Y_a}(S) \ar[ddl]_{\sigma} \ar[ddr]^{p} &  \\
 & & \\
X_a^0\subset \overline{X_a}\subset \PP \JA_0 \ar@{-->}[rr]^{\phi} &  &
Y_a\subset \PP U_a,}$$ 
where $p$ is the $\PP^1$-bundle defined by the restriction to $Y_a$ 
of the tautological rank two bundle $S$ over $G(2,\JA_0)$, 
and $\sigma$ is the blow-up of $\overline{X_a}$. 
\end{prop} 

\proof It follows from Proposition \ref{subspaces} that the projection $\s$ to $\PP \JA_0$ is an 
isomorphism over the complement of $\overline{X_a}$. Moreover, the fiber of a
point of  $\overline{X_a}$ is a $\PP^a$, which is mapped isomorphically by $p$ to a maximal
linear subspace of $Y_a$. This implies in particular that $\s^{-1}(\overline{X_a})$
is a smooth irreducible divisor $E$ in $Z_a$, which is itself smooth since $Y_a$ is smooth. 
By \cite{einsb}, Theorem 1, this is enough to ensure that $\s$ is the blow-up of 
$\overline{X_a}$. \qed 

\begin{coro} The variety of reductions $Y_a$ is a rational Fano manifold of index $a+1$,
with Picard group $Pic(Y_a)=\ZZ\cO(1)$. 
\end{coro}

\proof The claim on the Picard group is clear. To compute the index of $Y_a$, 
let again $E$ denote the exceptional divisor of $\sigma$, and $H$ the pull-back
of the hyperplane class. Since there are three $\PP^a$'s through the 
general point of $Y_a$, we have $E.f=3$ if $f$ denotes the class
of a fiber of $p$. Also $H.f=1$, hence $p^*\cO(1)=3H-E$.  
On the other hand it is easy to see that $\cO_S(1)\ot p^*\cO(1)=
H$, hence $\cO_S(1)=E-2H$. Now one can compute the canonical 
divisor of $Z_a$ in two ways:
$$\begin{array}{rcl}
K_{Z_a} & = & \sigma^*K_{\PP \JA_0}+aE=-(3a+2)H+aE \\
 & = & p^*(K_{Y_a}\otimes\det S)\otimes \cO_S(-2)=p^*K_{Y_a}+H-E,
\end{array}$$
so that $p^*K_{Y_a}=-(a+1)p^*\cO(1)$ and $K_{Y_a}=\cO(-a-1)$. Since the Picard group 
of $Y_a$ is generated by $\cO(1)$, this implies that $Y_a$ is a Fano manifold of index
$a+1$. 

The fact that it is rational follows from the diagram above: if $L$ is a hyperplane 
in $\PP \JA_0$, it is birational to its strict transform by $p$, which is itself 
birational to $Y_a$ via $\s$, since a general line of $Y_a$ meets $L$. \qed

\medskip\noindent {\it Remark}. The diagram of Proposition \ref{blowup} 
can be interpreted in terms of the study we made in section 2. Indeed, we 
can complete it as follows:\medskip

$$\xymatrix@1{  & Z_a=\PP_{Y_a}(S) \ar[dl]_{\sigma} \ar[dr]^{p} &  \\
X_a^0\subset \overline{X_a}\subset \PP \JA_0 \ar@{-->}[rr]^{\phi} &  &
Y_a\subset \PP U_a, \\
 \cup & & \cup \\
\PP H^0(X^a_I,\cF^a_I)^0 \ar@{>}[rr]^{zero-set\;map} & & Y_a^0
}$$ 

\bigskip\noindent
Remember from {\bf 2.8} and {\bf 2.9} that we defined on $X^a_I$ a homogeneous vector bundle
$\cF^a_I$, whose space of global sections was isomorphic to $\JA_I=\JA_0$ (we take $w=I$ 
here, which is harmless). On the complement $\PP H^0(X^a_I,\cF^a_I)^0$  of the discriminant
hypersurface, the zero-locus of a projective section was a triality subvariety $Z^a_{\e}$,
and we proved (see Proposition 2.15) that the family of these triality varieties was 
parametrized by the open orbit $Y_a^0$ of the variety of reductions. The map $\phi$ is 
nothing but the zero-set map considered as a rational map. 

\medskip We now turn to a different direction. We use Proposition \ref{blowup}
to compute the Betti numbers of the $Y_a$'s. 

\begin{coro} $Y_a$ has pure cohomology (i.e. its Hodge numbers $h^{p,q}(Y_a)=0$ for
$p\ne q$), and for $a\ge 2$ its Betti numbers can be deduced from those of $X_a$ 
by the formula
$$b_{2p}(Y_a)=\frac{1+(-1)^p}{2}+\sum_{0\le 2j<a}b_{2p-4j-2}(X_a).$$
The topological Euler charateristic of $Y_a$ is $e(Y_a)=3\frac{a(a+2)}{2}+1$, 
again for $a\ge 2$.
\end{coro}

\proof A simple computation, using the formulas giving the Betti numbers of 
a blow-up that can be found in \cite{gh}, page 605. \qed

\medskip The Betti numbers of $X_a$ present a nice regular pattern at least for $a\ge 2$:
$$b_{2p}(X_a)=\Bigg \{ \begin{array}{l} 1 \quad {\rm for}\; 0\le p<\frac{a}{2}\; 
{\rm or}\; \frac{3a}{2}< p\le 2a, \\ 
2 \quad {\rm for}\; \frac{a}{2}\le p<a\; {\rm or}\; a< p\le \frac{3a}{2}, \\ 
3\quad {\rm for}\; p=a. \end{array}.$$
In particular $e(X_a)=3a+3$. From this fact and the recursive formula of the Corollary
we can easily deduce the explicit Betti numbers of our varieties of reductions:
$$\begin{array}{clc}
p & 0\;  1\;  2\;  3\;  4\;  5\;  6\;  7\;  8\;  9\;  10\,  11\,  12\,  13\,  
14\,  15\,  16\,  17\,  18\,  19\,  20\,  21\,  22\,  23\,  24  & \\
b_{2p}(Y_1) & 1\;  1\;  1\;  1\;       &               \\
b_{2p}(Y_2) & 1\;  1\;  3\;  3\;  3\;  1\;  1\;   &              \\
b_{2p}(Y_4) & 1\;  1\;  2\;  3\;  4\; 5\; 5\; 5\; 4\; 3\;\; 2\;\; 1\,\;\; 1 &    \\
b_{2p}(Y_8) & 1\;  1\;  2\;  2\; 3\; 4\; 5\; 6\; 7\; 8\;\; 8\;\; 9\;\;\, 9\;\;\, 9\;\;\, 
8\;\;\, 8\;\; 7\;\;\, 6\;\;\, 5\;\;\, 4\;\; 3\;\;\, 2\;\;\, 2\;\; 1\;\;\, 1  & 
\end{array}$$
Correspondingly, the Euler characteristics are 
$$ e(Y_1)=4, \quad  e(Y_2)=13, \quad  e(Y_4)=37, \quad e(Y_8)=121.$$

\bigskip
Now we identify the rational map $\phi$. Since $p^*\cO(1)=3H-E$,
it must be defined by a system of cubics on $\PP \JA_0$. 

\begin{prop} 
The space of cubics on $\PP \JA_0$ vanishing on $\overline{X_a}$
is isomorphic to $U_a$.
\end{prop}

\proof A case-by-case verification with \cite{LiE} shows that the space of cubics 
on $\JA_0$ is 
$$S^3\JA_0^*=S^3\JA_0 = S^{(3)}\JA_0\oplus S^2\JA_0\oplus U_a.$$
The embedding of $S^2\JA_0$ inside $S^3\JA_0$ is given as follows:
to $A,B\in \JA_0$, we associate the cubic form $p_{A,B}(X)=
\trace(X(AX)(BX))$. One can check that these cubics vanish on $X_a^0$, as 
expected,  but not identically on $\overline{X_a}$.

The embedding of $U_a$ is deduced from the map 
$\wedge^2\JA_0\rightarrow S^3\JA_0$ defined as follows: to a
skew-symmetric form $\th$, we associate the cubic form $$p_{\th}(X)=
\th(X,X^2-\frac{1}{3}\trace(X^2)I).$$ Such cubics vanish on $\overline{X_a}$. 
Indeed, let $X$ be the projection of some $Z\in X_a$, that is $X=Z-\frac{z}{3}I$
and $Z^2=zZ$, where $z=\tr (Z)$. Then $X^2-\frac{1}{3} \tr(X^2)I=\frac{z}{3}X$, 
so clearly $p_{\th}(X)=0$. 

To conclude that the base locus of $U_a\subset S^3\JA_0\simeq S^3\JA_0^*$
is exactly $\overline{X_a}$, we first observe that these cubics cannot vanish 
identically on the discriminant hypersurface, which is irreducible of degree $6$. 
By Proposition \ref{J0orbits}, what remains to check is that they don't 
vanish identically on ${\rm Tan}^0(X_a^0)$. 
Remember from the proof of Proposition \ref{orbits} that the tangent space
of $X_a^0$ at 
$$\begin{array}{rcl}
Z & = &  \begin{pmatrix} 1 & i & 0 \\ i & -1 & 0 \\ 0 & 0 & 0 \end{pmatrix} \quad {\rm is} \\
T_ZX_a & = & \Bigg\{X=\begin{pmatrix} r-iu_0 & u & v \\ \bar{u} & r+iu_0 & iv \\ \bar{v} 
& i\bar{v} & 0 \end{pmatrix}, \quad r\in {\bf C}, u,v\in\AA \Bigg\}.
\end{array}$$
Now consider the line of diagonal matrices in $\JA_0$, seen as a point $\th\in U_a$. 
A straightforward computation shows that $p_{\th}(X)=-\frac{2}{3}iu_0q(Im(u))\ne 0$. 
\qed

\begin{coro}
The rational map $\phi$ is the map defined by the linear system  of cubics 
vanishing on the projected Severi variety $\overline{X_a}$.
\end{coro}

Finally, we use Proposition \ref{blowup} to  compute the degrees of the varieties of reductions. 
We also provide the degrees of the Grassmannians $G(2,3a+2)$, which are well-known to be
the Catalan numbers $\frac{1}{3a+1}\binom{6a}{3a}$, to show that although the degrees 
of the $Y_a$ can be quite big, they are relatively small compared to those of their 
ambient Grassmannians. 

\begin{theo}
The degrees of the varieties $Y_a$,  and of the Grassmannians $G(2,\JA_0)$ 
are:
$$\begin{array}{rclclcl}
\deg Y_1 & = & 5 & \hspace{1cm} &\deg G(2,5) & = & 5 \\
\deg Y_2 & = & 57 & \hspace{1cm} &\deg G(2,8) & = & 132 \\
\deg Y_4 & = & 12\,273 & \hspace{1cm} &\deg G(2,14) & = & 208\,012 \\
\deg Y_8 & = & 1\,047\,361\,761 & \hspace{1cm} & \deg G(2,26) & = & 1\,289\,904\,147\,324 
\end{array}$$\end{theo}

\proof Using Proposition \ref{blowup}, the degree of $Y_a$ can be computed once
the normal bundle of $\overline{X_a}$ in $\PP \JA_0$, 
and the Chow ring of $X_a$ are understood. The case of $a=2$ is explained in 
the next section. 

The most complicated case is of course that of $Y_8$. We give a detailed description 
of the Chow ring of the Cayley plane $X_8=\OO\PP^2$ in \cite{op2}, and show in 
this paper how the degree of $Y_8$ can be computed. \qed

\subsection{Varieties of reductions are compactifications of affine spaces}
The fact that $Y_1$ is a compactification of $\CC^3$ is due to Furushima \cite{furu},
who gave several geometric proofs of this property. We show that the varieties of reductions
are always compactifications of affine spaces (and indeed {\it minimal compactifications}, since
the Picard group is cyclic). Actually, this will directly follow from the 
fact that there is only a finite number of orbits. 

\begin{theo}\label{affine}
 The variety of reductions $Y_a$ is a compactification of $\CC^{3a}$.\end{theo}

\proof Since $Y_a$ is a smooth projective variety, it is enough to prove that $\SOA$ 
contains a one-dimensional torus $T$ acting on $Y_a$ with a finite number of fixed points. 
By the work of Byalinicki-Birula (\cite{bb}, Theorem 4.4), this will ensure that $Y_a$ contains 
a dense affine cell. More precisely, one can attach to each fixed point $z\in Y_a$ the subset of 
$Y_a$ defined as the union of the points that are attracted by $z$ through the action of 
$t\in T$, when $t$ tends to zero, and this provides a cell decomposition of the variety.  

For $a>1$, we have seen in Proposition \ref{s3} that $\SOA$ contains the triality group
$T(\AA)$, a reductive subgroup of 
maximal rank. We  choose a maximal torus $H$ of $SO(\AA)$ contained in $T(\AA)$. 

As a $\ft(\AA)$-module, $\JA_0=\CC^2\op\AA_1\op\AA_2\op\AA_3$. In particular 
the set of weights of $\JA_0$ is very easy to describe: there are $3a$ non-zero 
weights of multiplicity one (which are all conjugate under the action of the Weyl group), 
and the weight zero, whose multiplicity equals two. Let us denote by $L_1,\ldots,L_{3a}$
the one-dimensional weight spaces, and by $P$ the plane of weight zero. 

Let $T$ be a generic one-dimensional subtorus of $H$. Then the set $\PP \JA_0^T$ 
of fixed points of $T$ in $\PP \JA_0$ is the union of $3a$ points $e_1,\ldots,e_{3a}$,
and a projective line $d$. For the induced action on the Grassmannian $G(2,\JA_0)$, 
the set of fixed points $G(2,\JA_0)^T$ is then the union of $\binom{3a}{2}$ points, 
the lines $\overline{e_ie_j}$ in $\PP \JA_0$, another point, the line $d$ in 
$\PP \JA_0$, and $3a$ projective lines $d_i$, given by the set of lines in  
$\PP \JA_0$ joining $e_i$ to some point of $d$. 

We need to check that none of these lines is contained in $Y_a$. Suppose that $d_i\subset Y_a$. 
This would mean that $L\subset (\fsoa l_i)^{\perp}$. But $l_i$ is generated by a weight vector 
$X_i$ contained in some $\AA_j\subset \JA_0$: for example, if $j=3$, 
$$\begin{array}{rcl}
X_i & = & \begin{pmatrix} 0 & z & 0 \\ \bar{z} & 0 & 0 \\ 0 & 0 & 0 \end{pmatrix}, \quad {\rm so} \\
\quad \fsoa X_i & = & \Bigg\{
\begin{pmatrix} 2q(z,a_3) & g_3(z) & -za_1 \\ \overline{g_3(z)} & -2q(z,a_3) & -\bar{z}a_2 \\ 
\bar{a_1}\bar{z} & -\bar{a_2}z & 0 \end{pmatrix},\; g\in\ft(\AA),a_1,a_2,a_3\in\AA\Bigg\}.
\end{array}$$
Since the diagonal coefficients of these matrices are not identically zero, $\fsoa X_i$
is not orthogonal to $P$. 

Therefore the fixed lines $d_i$ are not contained in $Y_a$, which means that the action
of $T$ on $Y_a$ has a finite number of fixed points. \qed

\medskip We can analyze a little more carefully the set of fixed points of $T$ contained
in $Y_a$. First note that the line $d$ is certainly one of them. If $X_i$ is as in
the proof above, then $\fsoa X_i$ is orthogonal to exactly one point of $d$ (more
precisely, one of the three points of $d\cap\overline{X_a}$). This proves that 
$d_i$ cuts $Y_a$ at exactly one point, giving $3a$ fixed point of $T$ in $Y_a$.

Finally, we have to decide which of the lines $\overline{e_ie_j}$ are contained in 
$Y_a$. Suppose that $e_i$ is generated by the vector $X_i$,above, which we denote by 
$A_3(z)$. If $e_j$ also 
belong to $\AA_3$ and is represented by $X_j=A_3(z')$, then the line $\overline{e_ie_j}$
represents a point  of 
$Y_a$ if $q(z',g_3(z))=0$ for all $g\in\ft(\AA)$. But $g_3(z)$ can be any vector in 
$\AA$ orthogonal to $z$, so this would imply that $z$ and $z'$ are parallel, hence $e_i=e_j$. 
Now suppose that $e_j$ does not belong to $\AA_3$ like $e_i$, but for example to $\AA_2$. 
Then the line $\overline{e_ie_j}$ represents a point of $Y_a$ if $q(z',za_1)=0$ for all
$a_1\in\AA$. But $z$ is isotropic, so this means that $z'\in L(z)$. Since 
$L(z)$ is, like $z$, preserved by the torus action, it has a basis of eigenvectors of 
the torus, and we have $\l=\dim L(z)$ 
possibilities for the choice of $e_j$ in $\AA_2$, and also in $\AA_3$. This gives 
$3a\l$ new fixed points of $T$ in $Y_a$. Since $\l=0$ when $a=1$ and $\l=a/2$ when
$a\ge 2$, we get a total of $4$ fixed points when $a=1$, and $3a^2/2+3a+1$ fixed 
points when $a\ge 2$. Note that we know from  \cite{bb} that this number of fixed points 
is just the topological Euler characteristic of $Y_a$, which we have thus recovered. 

A more interesting consequence is the fact that $\CC^{3a}$ is the complement in $Y_a$ 
of a hyperplane section. 

\begin{prop} For $a>1$, let $x$ be a point of the closed orbit $Y_a^4$ of $Y_a$, and let $H_x$ 
denote the polar hyperplane. Then $Y_a-Y_a\cap H_x\simeq\CC^{3a}$. \end{prop}

\proof A case by case examination with the help of \cite{LiE} shows that the highest
weight of $U_a$ has exactly $3a^2/2$ conjugate under the Weyl group action: they are 
the weights $\m_i+\m_j$ of the lines $\overline{e_ie_j}$ we have just described. 
In particular, the highest weight of $U_a$ is of this type, and we can find a 
one-parameter subgroup whose attractive point in $\PP U_a$ is the line of highest 
weight. This point $x$, which belongs to the closed orbit in $Y_a$ (and can be 
chosen, by homogeneity, to be any point in $Y_a^4$), attracts the 
complement of a hyperplane of $\PP U_a$ generated by the remaining weight vectors.
This is precisely the polar hyperplane $H_x$ of $x$. If we restrict the action to $U_a$, 
the point $x$ attracts the complement to the hyperplane section $Y_a\cap H_x$, and 
by Bialiynicki-Birula's theorem this is a copy of $\CC^{3a}$.  \qed 

\medskip Note that a point $x$ of $Y_a^4$ is a special line $l_x$ on $X_a^0$. 
The polar subspace $L_x$ has codimension $2$  in $\PP \JA_0$, and we get:

\begin{prop} The singular section $Y_a\cap H_x$ is the Schubert variety of lines
that meet $L_x$, and belong to $Y_a$. \end{prop}

\medskip\noindent {\it Remark}. When $a=1$ the analysis is of course different, 
since the codimension $4$ orbit is empty. Nevertheless, one can check that 
Byalinicki-Birula's method still applies, giving a family of singular hyperplane
sections of $Y_1$ parametrized by the closed orbit (a sextic curve -- or equivalently
by the special lines on $Y_1$), whose 
complements are isomorphic to $\CC^3$. Nevertheless, Furushima proved that 
there exists another family of hyperplane sections of $Y_1$ (parametrized by the general 
lines), whose complements are
affine cells. We do not know whether a similar
phenomenon holds for the other varieties of reductions $Y_a$, $a>1$.

\section{The case of $\PP^2\times\PP^2$}

The case $a=2$, $\AA=\CC$ deserves special attention because the variety of reductions 
$Y_2$ is a smooth compactification of the space of independent
triples in $\PP^2$. But we already know such a compactification : 
the Hilbert scheme ${\rm Hilb^3}\PP^2$. 

\subsection{$Y_2$ and the Hilbert scheme}

The second Severi variety $X_2=\PP^2\times\PP^2$ is embedded into $\PP\JC=\PP M_3(\CC)$,
and is homogeneous under the action of $SL_3(\CC)=SL_3\times SL_3$. Its hyperplane section 
$X_2^0=\PP T_{\PP^2}=\FF(1,2;3)\subset \PP\JC_0=\PP \fsl_3$, the variety of complete flags 
in $\CC^3$, is homogeneous under $SO_3(\CC)=SL_3$. We thus have a coincidence between $\JC_0$
and $\fso_3(\CC)$, and the map $\wedge^2\JC_0=\wedge^2\fsl_3\lra \fso_3(\CC)=\fsl_3$ is 
just the Lie bracket. Therefore:

\begin{prop} The second variety of reductions $Y_2\subset G(2,\fsl_3)$, is the variety 
of abelian planes in $\fsl_3$. \end{prop}

The open $SL_3$-orbit of $Y_2$ is the set of planes of matrices which are diagonal 
in a given basis. This basis is unique up to multiplication by scalars, so that a 
point in $Y_2^0$ is actually just an  unordered triple of independent points in 
$\PP^2$. In particular, $Y_2$ is birational to ${\rm Hilb^3}\PP^2$, the punctual 
Hilbert scheme of length three subschemes of $\PP^2$. But we can be much more
precise:
  
\begin{theo}
The Hilbert scheme ${\rm Hilb}^3\PP^2$ has two extremal contractions. 
One is the Hilbert-Chow morphism onto $Sym^3\PP^2$. The variety 
$Y_2$ is isomorphic to the image of the other one. 
\end{theo}

\proof The extremal contraction of ${\rm Hilb}^3\PP^2$ which is not the 
Hilbert-Chow morphism was constructed in \cite{wlz} as follows: this
is the morphism 
$$\varphi_1 : {\rm Hilb}^3\PP^2 \lra G(3,S^2\CC^3)$$
mapping a length three subscheme $Z$ of $\PP^2$ to the dimension $3$ 
system  of conics that contains it. 

\begin{lemm}
$U_2\simeq \wedge^3(S^2\CC^3)$.
\end{lemm}

\proof The space $\wedge^3(S^2\CC^3)$ is generated by decomposable 
tensors of the form $e^2\wedge f^2\wedge g^2$. If we identify the 
dual of $\CC^3$ with its second wedge power, we can associate to 
such a tensor the $3$-dimensional subspace of $\fgl_3$ generated
by $e\wedge f\otimes g$, $f\wedge g\otimes e$ and  $g\wedge e\otimes f$. 
These three morphism commute and their sum is $e\wedge f\wedge g$
times the identity, hence their projection to $\fsl_3$ from the 
identity defines a point in $U_2$. The morphism so defined is 
clearly $\fsl_3$ and $H_2$ equivariant (see {\bf 3.1} for the definition of $H_2$), 
hence an isomorphism. \qed
 
\smallskip
This implies that we can identify $Y_2$ with the subvariety 
of $G(3,S^2\CC^3)$ defined as the closure of the planes 
$e^2\wedge f^2\wedge g^2$. This is not exactly the image of $\varphi_1$, 
which is the closure of the planes $ef\wedge fg\wedge ge$. 
But the theorem follows from the following lemma. 

\begin{lemm}
The endomorphism of $\wedge^3(S^2\CC^3)$ mapping 
$e^2\wedge f^2\wedge g^2$ to $ef\wedge fg\wedge ge$,
is an isomorphism. 
\end{lemm}

\smallskip\noindent {\it Proof of the lemma}. Let $\mu$ denote this endomorphism. 
We prove that $2\mu$ is an involution. We have $64\mu (ef\wedge fg\wedge ge)$ is equal to
$$\begin{array}{l}
 \mu \Big( [(e+f)^2-(e-f)^2]\wedge
[(f+g)^2-(f-g)^2]\wedge[(g+e)^2-(g-e)^2] \Big) \\
 \qquad = \;  \sum_{\e,\e',\e''=\pm 1}
(f+\e e)(g+\e' f)\wedge (g+\e' f)(e+\e'' g)\wedge (e+\e'' g)(f+\e e) \\
 \qquad = \; \sum_{\e,\e',\e''=\pm 1}(fg\wedge ge\wedge ef + ef\wedge fg\wedge ge) \\ 
\qquad = \;
16ef\wedge fg\wedge ge.
\end{array}$$
Indeed, in the previous sum we need only keep terms which have even degree in 
each of the $\e, \e', \e''$ (the other ones clearly cancel), and there are only 
two of them. \qed 

\medskip Note that we have obtained, by the way, another interpretation of $Y_2$. 

\begin{prop} The second variety of reductions $Y_2$ is isomorphic to the 
variety of trisecant planes to the Veronese surface $v_2(\PP^2)\subset\PP^5$. 
\end{prop}

\medskip
It was proved in \cite{wlz} that $\varphi_1$ contracts those subschemes
of $\PP^2$ contained in a line, to a $\PP^2$ parametrizing precisely 
these lines. At the level of $SL_3$-orbits, the Hilbert scheme contains 
seven orbits, and $3$ of them are contracted
onto one the two $\fsl_3$-invariant $\PP^2$'s inside $Y_2$. Actually,  
$\varphi_1$ is just the blow-up of this $\PP^2$. The fact that the image of $\varphi_1$ 
is smooth, which is the most surprising point here, is not mentionned 
in \cite{wlz}. 

\begin{prop}\label{57}
The degree of $Y_2$ is $57$.
\end{prop}

\noindent {\sl First proof}. From the preceding theorem we get that 
$$\deg Y_2=(\varphi_1^*\cO(1))^6,$$ and we are reduced to a computation on 
${\rm Hilb}^3\PP^2$, whose Chow ring has been described in detail in 
\cite{el}. With their notations, it is easy to see that $\varphi_1^*\cO(1)
=A+H$. Since they computed that $H^6=15$, $H^5A=15$, $H^4A^2=3$, 
$H^3A^3=-12$, $H^2A^4=12$, $HA^5=-3$ and $A^6=-15$, we get
$$\varphi_1^*\cO(1)^6  =  H^6+6H^5A+15H^4A^2+
20H^3A^3+15H^2A^4+6HA^5+A^6 =57.$$

\smallskip\noindent  {\sl Second proof}. We can use Theorem \ref{blowup} 
and the structure of  $Z_2$
to compute the degree of $Y_2$ as follows. First notice that $H.f=1$
implies that $$\deg Y_2=H(3H-E)^6$$ (recall that $H$ is the pull-back of 
the hyperplane class by $\sigma$, and $E$ the exceptional divisor). 
The intersection numbers on $Z_2$ can be computed explicitely, using 
the fact that $\sigma$ is a blow-up with smooth center $\overline{X_2}$, 
once we know the Chern classes of its normal bundle. First note that 
the normal bundle of $X_2=\PP^2\times\PP^2$ is 
$Q(1)\otimes Q'(1)$, where $Q, Q'$ denote 
the rank two tautological quotient bundles on the two copies of $\PP^2$.
Since $\overline{X_2}$ is an isomorphic linear projection of $X_2$, we have 
an exact sequence $0\ra \cO(1,1)\ra N_{X_2}\ra N_{\overline{X_2}}\ra 0$, 
from which we deduce that the Chern classes of $N_{\overline{X_2}}$ are
$$\begin{array}{rcl}
c_1(N_{\overline{X_2}}) & = & 5h+5h', \\
c_2(N_{\overline{X_2}}) & = & 10h^2+17hh'+10(h')^2, \\
c_3(N_{\overline{X_2}}) & = & 18hh'(h+h'), 
\end{array}$$
if $h$ and $h'$ denote the two hyperplane classes on our two copies of $\PP^2$. 
Therefore, the Chow ring of the exceptional divisor $E$ is the 
quotient of $\ZZ[h,h',e]$ by the relations $h^3=(h')^3=0$ and 
$e^3-5(h+h')e^2+(10h^2+17hh'+10(h')^2)e-18hh'(h+h')$ (see e.g. \cite{gh}). 
This being given, we compute (note that $h^2(h')^2e^2=1$, since $h^2(h')^2$ 
is the class of a fiber of $E\ra\PP^2\times\PP^2$, over which $e$ restricts
to the hyperplane class on $\PP^2$),
$$\begin{array}{l}
H^7=1, \quad H^6E=0, \quad H^5E^2=0, \quad H^4E^3=\deg\overline{X_2}=6,\\
H^3E^4=(h+h')^3e^3=5(h+h')^4e^2=30,\\
H^2E^5=(h+h')^2e^4=5(h+h')^3e^3-(h+h')^2(10h^2+17hh'+10(h')^2)e^2 \\
 \hspace{9cm} =150-54=96,\\
HE^6=(h+h')e^5=5(h+h')^2e^4-(h+h')(10h^2+17hh'+10(h')^2)e^3 \\
\hspace*{3.3cm} +18(h+h')(h^2h'+h(h')^2)e^2=480-270+36=246,
\end{array}$$
hence $\deg Y_2 = 3^6-3^3\binom{6}{3}6+3^2\binom{6}{2}30-
3\binom{6}{1}96+246=57$.  \qed

\begin{coro}
As a subvariety of $G(3,S^2\CC^3)$, the homology class of $Y_2$ is 
Poincar\'e dual to $\sigma_3+2\sigma_{21}+4\sigma_{111}$.
\end{coro}

\proof The Schubert class $\sigma_{111}$ is dual to the space of $3$-planes
in $S^2\CC^3$ contained in a fixed $4$-plane. This defines a $\PP^3$ 
cutting the Veronese surface in $\PP^5=\PP(S^2\CC^3)$ in four points. 
There are four ways to choose three of these four points, so that 
$\sigma_{111}$ cuts $Y_2$ at four points. 
 
Similarly, the Schubert class $\sigma_3$ is dual to the space of 
$3$-planes in $S^2\CC^3$ containing a fixed $2$-plane. For such a generic 
plane, there is a unique basis of $\CC^3$ diagonalizing each of the 
quadratic forms it parametrizes, so that $\sigma_2$ cuts $Y_2$ at 
a single point. 

To compute the last coefficient, we just note that a variety whose 
homology class is Poincar\'e dual to $x\sigma_3+y\sigma_{21}+z\sigma_{111}$
has degree $d=5x+16y+5z$. For $d=57$, $x=1$ and $z=4$ imply $y=2$. \qed

\subsection{A Calabi-Yau linear section}

Since $Y_2$ is a Fano variety of dimension $6$ and index $3$,
we can take smooth linear sections of dimension $3$ to get a family
of Calabi-Yau manifolds $C$. It was communicated to us by K. Ranestad that
these Calabi-Yau's had first been considered in \cite{tjotta}, who computed 
the corresponding Gromov-Witten invariants. The following result seems to be new:

\begin{prop}
The Betti numbers of $C$ are $b_0=b_2=b_4=b_6=1$, $b_1=b_5=0$, $b_3=2140$.
\end{prop}

\proof We know the Betti numbers of $Y_2$, and then all the Betti numbers of $C$
except $b_3$ are given by Lefschetz theorem. 
To compute $b_3$, we need to compute the Euler number
$$e(C)=\int_Cc_3(C)=\int_{Y_2}\frac{h^3}{(1+h)^3}c(Y_2),$$
where $h$ denotes the hyperplane class. To do this, we can pull everything 
back to ${\rm Hilb}^3\PP^2$. Note that since the the Hilbert scheme is the 
blow-up of $Y_2$ along a $\PP^2$ which will not be cut by a generic 
linear section of codimension $3$, we have 
$$e(C)=e(\varphi_1^{-1}(C))  =  \int_{{\rm Hilb}^3\PP^2}\frac{l^3}{(1+l)^3}
c({\rm Hilb}^3\PP^2),$$
where $l=\varphi_1^{-1}h$. To compute this, we use Bott's fixed-point 
formula as in \cite{es}, taking profit of the natural action 
on $\PP^2$, hence on ${\rm Hilb}^3\PP^2$, of the diagonal torus 
$D=\{diag(x_0,x_1,x_2)\}$ of $GL_3$. 

The fixed points in ${\rm Hilb}^3\PP^2$
are unions of monomial ideals supported on the three fixed points in $\PP^2$. 
There are $22$ of them, divided into five classes of cardinality 
$1$, $6$, $6$, $6$ and $3$ respectively. These classes are described below
with the induced action of the torus on the tangent space of the Hilbert
scheme at the fixed points. 

The drawings on the left column of the table below represent the
different types of length three subschemes of $\PP^2$ which are 
fixed points of the torus action: 1) the union of the three
fixed points, 2) a length two subscheme supported on a fixed point, 
given by a tangent line pointing to another fixed point, plus that 
point, 3) a length two subscheme supported on a fixed point, 
given by a tangent line pointing to another fixed point, plus the 
other fixed point, 4) a length three curvilinear subscheme supported
on a fixed point, with a prefered direction pointing to another 
fixed point, 5) a fattened fixed point, i.e. the square of the maximal 
ideal of a fixed point. 

For each fixed point $Z$, the torus $D$ acts on the tangent space $T_Z{\rm Hilb}^3\PP^2$. 
The character ${\rm ch}_D (T_Z{\rm Hilb}^3\PP^2)$ of this module has been computed with the 
help of formula (4.7) in \cite{es}. 

$$\begin{array}{ccc}
Z\quad\qquad & H^0({\mathcal O}_Z) & {\rm ch}_D (T_Z{\rm Hilb}^3\PP^2) \\
\setlength{\unitlength}{3mm}
\begin{picture}(8,5)(4,2)
\put(4,0){$\bullet$} \put(8,0){$\bullet$} \put(6,3.4){$\bullet$} 
\end{picture} & x_0^3, x_1^3, x_2^3 &
\frac{x_0}{x_1}+\frac{x_0}{x_2}+\frac{x_1}{x_0}+
\frac{x_1}{x_2}+\frac{x_2}{x_0}+\frac{x_2}{x_1} \\
\setlength{\unitlength}{3mm}
\begin{picture}(8,5)(4,2)
\put(6,3.4){$\circ$} \put(8,0){$\bullet$} \put(4,0){$\bullet$} 
\put(4.3,.3){\vector(1,0){2}}
\end{picture} & x_0^3, x_2^3, x_2^2x_0 &
2\frac{x_0}{x_1}+\frac{x_0}{x_2}+\frac{x_2}{x_0}+
\frac{x_2}{x_1}+(\frac{x_2}{x_0})^2 \\
\setlength{\unitlength}{3mm}
\begin{picture}(8,5)(4,2)
\put(6,3.4){$\bullet$} \put(8,0){$\circ$} \put(4,0){$\bullet$} 
\put(4.3,.3){\vector(1,0){2}}
\end{picture} & x_0^3, x_2^3, x_2^2x_1 &
\frac{x_0}{x_1}+\frac{x_1}{x_0}+\frac{x_0}{x_2}+
\frac{x_2}{x_0}+\frac{x_2}{x_1}+(\frac{x_2}{x_1})^2 \\
\setlength{\unitlength}{3mm}
\begin{picture}(8,5)(4,2)
\put(6,3.4){$\circ$} \put(8,0){$\circ$} \put(4,0){$\bullet$} 
\put(4.3,.4){\vector(1,0){2}}\put(4.3,.2){\vector(1,0){2}}
\end{picture} & x_0^3, x_0^2x_1, x_0x_1^2 &
\frac{x_0}{x_1}+\frac{x_1}{x_2}+\frac{x_0}{x_2}+
\frac{x_1^2}{x_0x_2}+(\frac{x_0}{x_1})^2+(\frac{x_0}{x_1})^3 \\
\setlength{\unitlength}{3mm}
\begin{picture}(8,5)(4,2)
\put(6,3.4){$\circ$} \put(8,0){$\circ$} \put(4,0){$\bullet$} 
\put(4.3,.3){\vector(1,0){2}}\put(4.25,.4){\vector(2,3){1}}
\end{picture} & x_0^3, x_0^2x_1, x_0^2x_2 &
2\frac{x_0}{x_1}+\frac{x_0x_1}{x_2^2}+2\frac{x_0}{x_2}+
\frac{x_0x_2}{x_1^2}
\end{array}$$

\bigskip

\bigskip\noindent 
\smallskip We now choose a one dimensional 
subtorus 
$$T=\Bigg\{ \begin{pmatrix} t^{w_0} & 0 & 0 \\ 0 & t^{w_1} & 0 \\ 0 & 0 &  t^{w_2}
\end{pmatrix}, t\in {\bf C}^*\Bigg\}\subset 
D=\Bigg\{ \begin{pmatrix} x_0 & 0 & 0 \\ 0 & x_1 & 0 \\ 0 & 0 &  x_2
\end{pmatrix}, x_0,x_1,x_2\in {\bf C}^*\Bigg\},$$ 
with the same fixed points. 
This condition is achieved when the tangent spaces have no zero 
weight space: the formulas above show that the weights $w_0, w_1, w_2$ 
must be such that $w_i\neq w_j$ and $2w_i\neq w_j+w_k$ for all
permutations $i,j,k$ of $0,1,2$. We let $w_0=0, w_1=1, w_2=3$. 
The table below gives the character ${\rm ch}_T (T_Z{\rm Hilb}^3\PP^2)$,
that is, the set of integers $m_1,\ldots ,m_6$ such that
$${\rm ch}_T (T_Z{\rm Hilb}^3\PP^2)=t^{m_1}+t^{m_2}+t^{m_3}+t^{m_4}+t^{m_5}+t^{m_6}.$$

To complete the computation, we will apply Bott's fixed point formula in the 
form given by Theorem 2.2 in \cite{es}. The intersection numbers we need to compute are 
given by this formula as sums of rational numbers which are contributions of the 
$22$ fixed points. Part of these contributions come from the Chern classes of 
the tangent bundle, and are easily deduced from its character at each fixed point. 
The contribution of the class $l$ can be computed as follows: recall that $l=\varphi_1^{-1}h
=2H+E$, where $H$ denotes the pull-back of the hyperplane class by the Hilbert-Chow morphism, 
and $E$ is the exceptional divisor of that morphism. If $\cE_n$ is the vector bundle whose 
fiber at a scheme $Z$ is given by $H^0(\cO_Z(n))$, we have $\det\cE_n=
nH+E$, so we just need to compute $\det\cE_0$ and $\det\cE_1$, which is straightforward. 
We get: 

$$\begin{array}{cccc}
Z\qquad & \quad\det\cE_0\quad  & \quad \det\cE_1\quad  &\quad  l\quad  \\
\setlength{\unitlength}{2mm}
\begin{picture}(8,5)(4,2)
\put(4,0){$\bullet$} \put(8,0){$\bullet$} \put(6,3.4){$\bullet$} 
\end{picture} & 1 & x_0x_1x_2 & x_0^2x_1^2x_2^2 \\ 
\setlength{\unitlength}{2mm}
\begin{picture}(8,5)(4,2)
\put(6,3.4){$\circ$} \put(8,0){$\bullet$} \put(4,0){$\bullet$} 
\put(4.3,.3){\vector(1,0){2}}
\end{picture} & \frac{x_0}{x_2} & x_0^2x_2 & x_0^3x_2^3 \\
\setlength{\unitlength}{2mm}
\begin{picture}(8,5)(4,2)
\put(6,3.4){$\bullet$} \put(8,0){$\circ$} \put(4,0){$\bullet$} 
\put(4.3,.3){\vector(1,0){2}}
\end{picture} & \frac{x_2}{x_1} & x_0x_1x_2 & x_0^2x_1^3x_2 \\
\setlength{\unitlength}{2mm}
\begin{picture}(8,5)(4,2)
\put(6,3.4){$\circ$} \put(8,0){$\circ$} \put(4,0){$\bullet$} 
\put(4.3,.4){\vector(1,0){2}}\put(4.3,.2){\vector(1,0){2}}
\end{picture} & \frac{x_1^3}{x_0^3} & x_1^3 & x_0^3x_1^3 \\
\setlength{\unitlength}{2mm}
\begin{picture}(8,5)(4,2)
\put(6,3.4){$\circ$} \put(8,0){$\circ$} \put(4,0){$\bullet$} 
\put(4.3,.3){\vector(1,0){2}}\put(4.25,.4){\vector(2,3){1}}
\end{picture} & \frac{x_1x_2}{x_0^2} & x_0x_1x_2 & x_0^4x_1x_2.
\end{array}$$

\medskip
Finally, we list in the table below, 
for our choice of $T$, and for each fixed point $Z$, the integers 
$m_1,  \ldots ,m_6$ giving the character of the tangent space, the 
integers $c_i$ corresponding to the Chern classes $c_i({\rm Hilb}^3\PP^2)$ for 
$i=1,2,3,6$ (these are just the $i$-th elementary symmetric functions of $m_1,\ldots ,m_6$), 
as well as the weights $\l$ of $l$. 

From these data we can compute that 
$$\int_{{\rm Hilb}^3\PP^2}c_3({\rm Hilb}^3\PP^2)l^3 = \sum_Z\frac{c_3}{c_6}\l^3 = 243,$$

$$\int_{{\rm Hilb}^3\PP^2}c_2({\rm Hilb}^3\PP^2)l^4 = \sum_Z\frac{c_2}{c_6}\l^4 = 261,$$

$$\int_{{\rm Hilb}^3\PP^2}c_1({\rm Hilb}^3\PP^2)l^5 = \sum_Z\frac{c_1}{c_6}\l^5 = -171,$$

\noindent hence $e(C)=243-3\times 261-6\times 171-10\times 57=-2136$. Since 
$e(C)=4-b_3$, this implies our claim. \qed
 
\medskip
Since $b_3(C)=2+2h^{2,1}(C)$, we get $h^{2,1}(C)=1069=\dim H^1(C,TC)$. This is 
the dimension of the space of deformations of $C$, and is much larger than the number
of available parameters for the codimension three linear section which defines $C$. 
Therefore:

\begin{coro} A general deformation of $C$ is not a linear section of $Y_2$. \end{coro}

{\small 
$$\begin{array}{crrrrr}
T_Z{\rm Hilb}^3\PP^2 & c_1 & c_2 & c_3 & c_6 & \l \\
 & &  & & & \\
1,-1,2,-2,3,-3 & 0 & -14 & 0 & -36 & 8 \\ 
 & &  & & & \\
-1, -1, 2, 3, -3, 6 & 6 & -12 & -70 & -108 & 9 \\
1,1,2,-2, 3, 4 & 9 & 23 & -5 & -48 & 12 \\
1, -1, 2, -2, -3, -3 & -6 & 4& 30 & 36 & 3 \\
1, -1, -2, -2, -2, -3 & -9 & 29 & -35 & -24 & 3 \\
1, 2, -2, 3, 3, -4 & 3 & -17 & -63 & 144 & 12 \\
-1, 2, 2, 3, -3, -6 & -3 & -27 & 23 & -216 & 9\\
 & &  & & & \\
1, -1, 2, 3, -3, 4 & 6 & -2& -60 & 72 & 10\\
1, -1, -2, 3, -3, -4 & -6 & -2& 60 & 72 & 6\\
1, -1, 2, -2, 3, 6 & 9 & 13& -45 & 72 & 11\\
1, -1, 2, -2, -3, -6 & -9 & 13& 45 & 72 & 5\\
1, 2, 2, -2, 3, -3 & 3 & -11& -39 & 72 & 9 \\
-1, 2, -2, -2, 3, -3 & -3 & -11& 39 & 72 & 7 \\
 & &  & & & \\
-1, -1, -2, -2, -3, -3 & -12 & 58 & -144 & 36 & 3 \\
1, 2, -2, 3, -3, -4 & -3 & -17& 39 & -144 & 3 \\
-1, 2, -3, 5, -6, -9 & -12 & -6& 368 & 1620 & 9 \\
-1, -2, 3, -4, 5, -6 & -3 & -41& 87 & -720 & 12 \\
-1, 2, 3, -4, 6, 9 & 15 & 39 & -235 & 1296 & 9 \\
1, -1, 2, 3, 4, 6 & 15 & 79 & 165 & -144 & 12 \\
 & &  & & & \\
1, -1, -1, -3, -3, -5 & -12 & 49 & -72 & -45 & 4 \\
1, 1, -2, -2, 4, -5 & -3 & -21 & 47 & -80 & 7 \\
1, 2, 2, 3, 3, 4 & 15 & 91 & 285 & 144 & 13 
\end{array}$$
}

\bigskip\noindent 
Atanas ILIEV, Institut de Math\'ematiques, Acad\'emie des Sciences 
de Bulgarie, 8 Rue Acad. G. Bonchev, 1113 Sofia, Bulgarie.

\noindent Email : ailiev@math.bas.bg

\smallskip\noindent 
Laurent MANIVEL, Institut Fourier, Laboratoire de Math\'ematiques, 
UMR 5582 (UJF-CNRS), BP 74, 38402 St Martin d'H\`eres Cedex, France.

\noindent Email : Laurent.Manivel@ujf-grenoble.fr

\end{document}